\documentstyle[amscd]{amsart}

\def\N{{\Bbb N}}
\def\Z{{\Bbb Z}}
\def\Q{{\Bbb Q}}
\def\R{{\Bbb R}}
\def\C{{\Bbb C}}

\def\V{{\Bbb V}}
\def\E{{\Bbb E}}
\def\B{{\Bbb B}}
\def\bH{{\Bbb H}}
\def\K{{\Bbb K}}
\def\bL{{\Bbb L}}
\def\bS{{\Bbb S}}
\def\bT{{\Bbb T}}

\def\A{{\cal A}}
\def\H{{\cal H}}
\def\I{{\cal I}}
\def\J{{\cal J}}
\def\M{{\cal M}}
\def\T{{\cal T}}
\def\L{{\cal L}}
\def\calB{{\cal B}}
\def\calC{{\cal C}}
\def\calK{{\cal K}}
\def\calV{{\cal V}}
\def\calX{{\cal X}}
\def\calZ{{\cal Z}}

\def\X{\overline{X}}
\def\Sbar{\overline{S}}
\def\taut{\tilde{\tau}}
\def\tauh{\hat{\tau}}
\def\w{\omega}
\def\wv{\vec{w}}
\def\f{\vec{f}}

\def\h{\goth h}
\def\del{\partial}
\def\delbar{\overline{\del}}
\def\etilde{\tilde{e}}

\def\itilde{\tilde{\imath}}
\def\stilde{\tilde{s}}
\def\pbar{\overline{p}}

\def\dot{\bullet}
\def\blank{\underline{\phantom{x}}}

\def\com{[\phantom{x},\phantom{x}]}
\def\id{\text{id}}
\def\sgn{\text{sgn}}

\def\rel#1{{\calC_g/#1}}

\newsymbol \ltimes 226E 
\newsymbol \rtimes 226F
\newsymbol \boxplus 1201

\newcommand\coker{\operatorname{coker}}
\newcommand\Hom{\operatorname{Hom}}
\newcommand\Ext{\operatorname{Ext}}
\newcommand\Aut{\operatorname{\text{Aut}}}
\newcommand\Diff{\operatorname{Diff}}
\newcommand\Jac{\operatorname{Jac}}
\newcommand\Pic{\operatorname{Pic}}

\newtheorem{theorem}{Theorem}[section]
\newtheorem{lemma}[theorem]{Lemma}
\newtheorem{proposition}[theorem]{Proposition}
\newtheorem{corollary}[theorem]{Corollary}

\theoremstyle{definition}

\newtheorem{definition}[theorem]{Definition}

\theoremstyle{remark}

\newtheorem{remark}[theorem]{Remark}

\begin{document}

\title[The Torelli Groups and  Geometry]{Torelli Groups and
Geometry of Moduli Spaces of Curves}

\author{Richard M.~Hain}

\address{Department of Mathematics\\
Duke University\\ Durham, NC 27708-0320}

\thanks{Research  supported in part by grants from the National
Science Foundation and the NWO.}

\email{hain@@math.duke.edu}

\maketitle

\section{Introduction}\label{intro}

The Torelli group $T_g$ is the kernel of the natural homomorphism
$\Gamma_g \to Sp_g(\Z)$ from the mapping class group in genus $g$ to
the group of $2g \times 2g$ integral symplectic matrices. It accounts
for the difference between the topology of $\A_g$, the moduli space of
principally polarized abelian varieties of dimension $g$, and $\M_g$,
the moduli space of smooth projective curves of genus $g$, and
therefore should account for some of the difference between their
geometeries. For this
reason, it is an important problem to understand its structure and its
cohomology. To date, little is known about $T_g$ apart from Dennis
Johnson's few fundamental results --- he has proved that $T_g$ is
finitely generated when $g\ge 3$ and has computed $H_1(T_g,\Z)$. It
is this second result which will concern us in this paper.  Crudely
stated, it says that there is an $Sp_g(\Z)$-equivariant isomorphism
$$
H^1(T_g,\Q) \approx PH^3(\Jac C,\Q)
$$
where $C$ is a smooth projective curve of genus $g$, and $P$ denotes
primitive part. My aim in this paper is to give a detailed exposition
of Johnson's homomorphism
$$
PH^3(\Jac C,\Q) \to H^1(T_g,\Q)
$$
and to explain how Johnson's computation, alone and in concert
with M.~Saito's theory of Hodge modules \cite{saito}, has some
remarkable consequences
for the geometry of $\M_g$. It implies quite directly, for example, that
for each $l$, the Picard group of the moduli space $\M_g(l)$ of genus
$g\ge 3$ curves with a level $l$ structure is finitely generated. Combined
with Saito's work, it enables one to completely write down all ``natural''
generically defined normal functions over $\M_g(l)$ when $g\ge 3$.
The result is that modulo torsion, all are half integer multiples of the
normal function of the cycle $C-C^-$. This is applied to give a new
proof of the Harris-Pulte Theorem \cite{harris,pulte}, which relates
the mixed Hodge structure on the fundamental group of a curve $C$ to
the algebraic cycle $C-C^-$ in its jacobian. We are also able to
``compute'' the archimedean height pairing between any two ``natural''
cycles in a smooth projective variety defined over the moduli space of
curves, provided they are homologically trivial over each curve, disjoint
over the generic curve, and satisfy the usual dimension restrictions.
The precise statement can be found in Section \ref{heights}.

The classical Franchetta conjecture asserts that the Picard group
of the generic curve is isomorphic to $\Z$ and is generated by
the canonical divisor. Beauville (unpublished),  and later Arbarello
and Cornalba \cite{arb-corn}, deduced this from Harer's computation of
$H^2(\Gamma_g)$. As another application of the classification of
normal functions over $\M_g(l)$, we prove a ``Franchetta Conjecture''
for the generic curve with a level $l$ structure. The statement is that
the Picard group of the generic curve of genus $g$ with a level $l$
structure is finitely generated of rank 1 --- the torsion subgroup is
isomorphic to $(\Z/l\Z)^{2g}$; mod torsion, it is generated by the
canonical bundle if $l$ is odd, and by a square root of the canonical
bundle if $l$ is even. Our proof is only valid when $g\ge 3$;
it does not use the computation of $\Pic \M_g(l)$, which is
not known at this time. We also compute the Picard group of the generic
genus $g$ curve with a level $l$ structure and $n$ marked points.

Our results on normal functions are inspired by
those in the last section of Nori's remarkable paper \cite{nori} where
normal functions on finite covers of Zariski open subsets of
the moduli space of principally polarized abelian varieties are studied.
There are analogues of our main results for $\A_g(l)$, the moduli
space of principally polarized abelian varieties of dimension $g$
with a level $l$ structure. These results are similar to Nori's,
but differ. The detailed statements, as well as a discussion of the
relation between the results, are in Section \ref{abelian}. Our
results on abelian varieties seem to be related to some results of
Silverberg \cite{silverberg}.

Sections \ref{homom} and \ref{original} contain an exposition of the
three constructions of the Johnson homomorphism that are given in
\cite{johnson:survey}. Since no proof of their equivalence appears
in the literature, I have given a detailed exposition, especially
since the equivalence of two of these constructions is essential in
one of the applications to normal functions.

In Section \ref{root} Johnson's results is used to give an explicit
description of the quotient of $\Gamma_g$  by the kernel of its action
on all $n$th roots of the canonical bundle is given. A consequence
of this computation is the ``well known fact'' that the only roots of
the canonical bundle defined over Torelli space are the canonical bundle
itself and all
theta characteristics.
\medskip

\noindent{\it Acknowledgements.}\ First and foremost, I would like to
thank Eduard Looijenga for his hospitality and for stimulating discussions
during a visit to the University of Utrecht in the spring of 1992 during
which some of the work in this paper was done. I would also like to thank
the  University of Utrecht and the Dutch NWO for their generous support
during that visit. Thanks also to Pietro  Pirola and Enrico Arbarello for
pointing out to me that the non-existence of sections of the universal
jacobian implies the classical Franchetta conjecture. From this it was a
short step to the generalizations in  Section \ref{franchetta}.

\section{Mapping Class Groups and Moduli}

At this time there is no argument within algebraic geometry to
compute the Picard groups of all $\M_g$, and one has to resort to
topology to do this computation.  Let $S$ be a compact orientable
surface of genus $g$ with $r$ boundary components and let $P$ be an
ordered set of $n$ distinct marked points of $S -\partial S$. Denote
the group of orientation preserving diffeomorphisms of $S$ that fix
$P\cup \partial S$ pointwise by $\Diff ^+ (S,P \cup \partial S)$.
Endowed with the compact open topology, this is a topological group. The
mapping class group $\Gamma_{g,r}^n$ is defined to be its group of
path components:
$$
\Gamma_{g,r}^n = \pi_0 \Diff^+(S,P\cup\partial S).
$$
Equivalently, it is the group of isotopy classes of orientation
preserving diffeomorphisms of $S$ that fix $P\cup \partial S$
pointwise.
It is conventional to omit the decorations $n$ and $r$ when they
are zero. So, for example, $\Gamma_g^n = \Gamma_{g,0}^n$.

The link between moduli spaces and mapping class groups is
provided by Teichm\"uller theory. Denote the moduli space of
smooth genus $g$ curves with $n$ marked points by $\M_g^n$.
Teichm\"uller theory provides a contractible complex manifold
$\calX_g^n$ on which $\Gamma_g^n$ acts properly discontinuously---it
is the space of all complete hyperbolic metrics on $S-P$ equivalent
under diffeomorphisms isotopic to the identity. The
quotient $\Gamma_g^n\backslash \calX_g^n$ is analytically
isomorphic to $\M_g^n$. It is useful to think of $\Gamma_g^n$ as
the orbifold fundamental group of $\M_g^n$.

One can compactify $S$ by filling in the $r$ boundary components of
$S$ by attaching disks. Denote the resulting genus $g$ surface by
$\Sbar$. Elements of $\Gamma_{g,r}^n$ extend canonically to $\Sbar$ to
give a homomorphism $\Gamma_{g,r}^n \to \Gamma_g$. Denote the
composite
$$
\Gamma_{g,r}^n \to \Gamma_g \to \Aut H_1(\Sbar,\Z)
$$
by $\rho$. Since elements of $\Gamma_{g,r}^n$ are represented by
orientation preserving diffeomorphisms, each element of
$\Gamma_{g,r}^n$ preserves the intersection pairing
$$
q : \Lambda^2 H_1(\Sbar,\Z) \to \Z.
$$
Consequently, we obtain a homomorphism
$$
\rho : \Gamma_{g,r}^n \to \Aut (H_1(\Sbar,\Z),q) \approx Sp_g(\Z).
$$
This homomorphism is well known to be surjective.

Denote the moduli space of principally polarized abelian varieties
of dimension $g$ by $\A_g$. Since this is the quotient of the Siegel
upper half plane by $Sp_g(\Z)$, it is an orbifold with orbifold
fundamental group $Sp_g(\Z)$.  The period map
$$
\M_g^n \to \A_g
$$
is a map of orbifolds and induces $\rho$ on fundamental groups.

The Torelli group $T_{g,r}^n$ is the kernel of the homomorphism
$$
\rho : \Gamma_{g,r}^n \to Sp_g(\Z).
$$
Since $\rho$ is surjective, we have an extension
$$
1 \to T_{g,r}^n \to \Gamma_{g,r}^n \to Sp_g(\Z) \to 1.
$$
The Torelli group $T_g$ encodes the differences between the topology of
$\M_g$ and $\A_g$ --- between curves and abelian varieties. More formally,
we have the Hochschild-Serre spectral sequence
$$
H^s(Sp_g(\Z),H^t(T_{g,r}^n)) \implies H^{s+t}(\Gamma_{g,r}^n).
$$

Much more (although not enough) is known about the topology of the
$\A_g$ than about that of the $\M_g$. For example, the rational
cohomology groups of the $\A_g$ stabilize as $g \to \infty$, and this
stable cohomology is known by Borel's work \cite{borel:triv}; it is
a polynomial ring generated by classes $c_1, c_3, c_5, \dots$, where
$c_k$ has degree $2k$. As with $\A_g$, the rational cohomology of
the $\M_g$ is known to stabilize, as was proved by Harer
\cite{harer:stab}, but the stable cohomology of the $\M_g$ is known
only up to dimension 4; the computations are due to Harer
\cite{harer:h3,harer:h4}.

Torelli space $\T_g^n$ is the quotient $T_g^n\backslash\calX_g^n$ of
Teichm\"uller space.  When $g\ge 3$, it is the moduli space of smooth
projective curves $C$, together with $n$ ordered distinct points and a
symplectic basis of $H_1(C,\Z)$.

The Torelli group is torsion free. Perhaps the simplest way to see this
is to note that, by standard topology, since $\calX_g^n$ is contractible,
each element of $\Gamma_g^n$
of prime order must fix a point of $\calX_g^n$. If $\phi\in \Gamma_g^n$
fixes the point corresponding to the marked curve $C$, then there is
an automorphism of $C$ that lies in the mapping class $\phi$. Since
the automorphism group of a compact Riemann surface injects into
$\Aut H^0(C,\Omega_C^1)$, and therefore into $H_1(C)$, it follows that
$T_g^n$ is torsion free. Because of this, the Torelli space
$\T_g^n$ is the classifying space of $T_g^n$.

One can view Siegel space $\h_g$ as the classifying space of
principally polarized abelian varieties of dimension $g$ together with
a symplectic basis of $H^1$. The period map therefore induces a
map
$$
\T_g^n \to \h_g
$$
which is 2:1 when $g\ge 2$, and ramified along the hyperelliptic
locus when $g\ge3$.

For a finite index subgroup $L$ of $Sp_g(\Z)$, let $\Gamma_{g,r}^n(L)$ be
the inverse image of $L$ in $\Gamma_{g,r}^n$ under the canonical
homomorphism $\Gamma_{g,r}^n \to Sp_g(\Z)$. It may be expressed
as an extension
$$
1 \to T_{g,r}^n \to \Gamma_{g,r}^n(L) \to L \to 1.
$$
Set $\M_{g}^n(L) = \Gamma_{g}^n(L)\backslash \calX_g^n$.
We will call $\Gamma_{g,r}^n(L)$ the {\it level $L$} subgroup of
$\Gamma_{g,r}^n$, and we will say that points in $\M_{g}^n(L)$
are curves with a level $L$ structure and $n$ marked points.
The traditional moduli space of curves with a level $l$ structure,
where $l\in \N^+$,  is obtained by taking $L$ to be the elements of
$Sp_g(\Z)$ that are congruent to the identity mod $l$.

Since the Torelli groups are torsion free, $\Gamma_{g,r}^n(L)$ is
torsion free when $L$ is. Note, however, that by the Lefschetz fixed
point formula, $\Gamma_{g,r}^n$ is torsion free when $n+2r > 2g+2$,
so that $\Gamma_{g,r}^n(L)$ may be torsion free even when $L$ is not.

\begin{proposition}
For all $g,n \ge 0$ and for each finite index subgroup $L$ of
$Sp_g(\Z)$, there is a natural homomorphism
$$
H^\dot(\M_g^n(L),\Z) \to H^\dot(\Gamma_g^n(L),\Z)
$$
which is an isomorphism when $\Gamma_g^n(L)$ is torsion free, and
is an isomorphism after tensoring with $\Q$ for all $l$.
\end{proposition}

\begin{pf}
Set $\Gamma=\Gamma_g^n$, $\Gamma(L) = \Gamma_g^n(L)$,
$\M(L)=\M_g^n(L)$ and $\calX=\calX_g^n$. Let $E\Gamma$ be any
space on which $\Gamma$ acts freely and properly discontinuously%
---so $E\Gamma$ is the universal covering space of some model of
the classifying space of $\Gamma$. Since $\calX$ is contractible, the
quotient $E\Gamma\times_{\Gamma(L)}\calX$ of $E\Gamma \times
\calX$ by the diagonal action of $\Gamma(L)$ is a model
$B\Gamma(L)$ of the classifying space of $\Gamma(L)$. The
projection $E\Gamma \times \calX \to \calX$ induces a map
$f: B\Gamma(L) \to \M(L)$ which
induces the map of the theorem. If $\Gamma(L)$ is torsion free, $f$ is a
homotopy equivalence. Otherwise, choose a finite index, torsion
free normal subgroup $L'$ of $L$. Then $\Gamma(L)$ is torsion free.
Set
$$
G = \Gamma(L)/\Gamma(L')\approx L/L'
$$
This is a finite group. We have the commutative diagram of Galois
$G$-coverings
$$
\begin{matrix}
B\Gamma(L') & \to & \M(L') \cr
\downarrow & & \downarrow \cr
B\Gamma(L) & \to & \M(L)
\end{matrix}
$$
where the top map is a homotopy equivalence. It therefore induces
a $G$-equivariant isomorphism
$$
H^\dot(\M(L')) \to H^\dot(B\Gamma(L')).
$$
The result follows as the vertical projections induce isomorphisms
$$
H^\dot(\M(L),\Q) \stackrel{\sim}{\to} H^\dot(\M(L'),\Q)^G
\text{ and }
H^\dot(\Gamma(L),\Q) \stackrel{\sim}{\to} H^\dot(\Gamma(L'),\Q)^G.
$$
\end{pf}

The group $\Gamma_{g,r}^n(L)$ also admits a moduli interpretation
when $r > 0$, even though algebraic curves have no boundary
components. The idea is that a topological boundary component of a
compact orientable surface should correspond to a first order local
holomorphic coordinate about a cusp of a smooth algebraic curve.
Denote by $\M_{g,r}^n(L)$ the moduli space of smooth curves of
genus $g$ with a level $L$ structure and with $n$ distinct marked
points and $r$ distinct, non-zero cotangent vectors, where the
cotangent vectors do not lie over any of the marked points, and
where no two of the cotangent vectors are anchored at the same
point. This is a $\left(\C^\ast\right)^r$ bundle over $\M_g^{r+n}(L)$.

\begin{proposition}
For all finite index subgroups $L$ of $Sp_g(\Z)$ and
for all $g,n,r \ge 0$, there is a natural homomorphism
$$
H^\dot(\M_{g,r}^n(L),\Z) \to H^\dot(\Gamma_{g,r}^n(L),\Z)
$$
which is an isomorphism when $\Gamma_{g,r}^n(L)$ is torsion free,
and is an isomorphism after tensoring with $\Q$ for all $L$.\qed
\end{proposition}

\section{The Johnson Homomorphism}
\label{homom}

Dennis Johnson, in a sequence of pioneering papers
\cite{johnson_1,johnson_2,johnson_3}, began a systematic
study of the Torelli groups. From the point of view of
computing the cohomology of the $\M_g$, the most important of
his results is his computation of $H_1(T_g^1)$ \cite{johnson_3}.
Let $S$ be a compact oriented surface of genus $g$ with a
distinguished base point $x_0$.

\begin{theorem}\label{john:h1}
There is an $Sp_g(\Z)$-equivariant homomorphism
$$
\tau_g^1 : H_1(T_g^1,\Z) \to \Lambda^3 H_1(S)
$$
which is an isomorphism mod 2-torsion.
\end{theorem}

Johnson has also computed $H_1(T_g^1,\Z/2\Z)$. It is related to
theta characteristics. Bert van~Geemen has interesting ideas regarding
its  relation to the geometry of curves.

A proof of Johnson's theorem is beyond the scope of this paper.
However, we will give three constructions of the homomorphism
$\tau_g^1$ and establish their equality.

We begin by sketching the first of these constructions: Since the Torelli
group is torsion free, there is a universal curve
$$
\calC \to \T_g^1
$$
over Torelli space. This has a tautological section
$\sigma: \T_g^1 \to \calC$.
There is also the jacobian
$$
\J \to \T_g^1
$$
of the universal curve. The universal curve can be imbedded in its
jacobian using the section $\sigma$---the restriction of this mapping
to the fiber over the point of Torelli space corresponding to $(C,x)$
is the Abel-Jacobi mapping
$$
\nu_x : (C,x) \to (\Jac C,0)
$$
associated to $(C,x)$. Since $T_g^1$ acts trivially on the first
homology of the curve, the local system associated to $H_1(C)$ is
framed. There is a corresponding topological trivialization of the
jacobian bundle:
$$
\J \stackrel{\sim}{\to} \T_g^1 \times \Jac C.
$$
Let $p : \J \to \Jac C$ be the corresponding projection onto the fiber.
Each element $\phi$ of $H_1(T_g^1,\Z)$ can be represented by an
imbedded circle $\phi :S^1 \to \T_g^1$. Regard the
universal curve
$\calC$ as  subvariety of $\J$ via the Abel-Jacobi mapping. Then the part
of the universal curve $M(\phi)$ lying over the circle $\phi$ is a 3-cycle
in $\J$. The Johnson homomorphism is defined by
$$
\tau_g^1(\phi) = p_\ast[M(\phi)] \in H_3(\Jac C,\Z) \approx
\Lambda^3 H_1(C,\Z).
$$

This definition is nice and conceptual, but is not so easy to work with.
In the remainder of this section, we remake this definition without
appealing to Torelli space. In the next section, we will give two more
constructions of it, both due to Johnson, and prove all three constructions
agree.

Recall that the {\it mapping torus} of a diffeomorphism $\phi$ of a
manifold $S$ is the quotient $M(\phi)$ of $S\times [0,1]$ obtained
by identifying $(x,1)$ with $(\phi(x),0)$:
$$
M(\phi) = S \times [0,1]\left\{(x,1) \sim (\phi(x),0)\right\}.
$$
The projection $S \times [0,1] \to [0,1]$ induces a bundle projection
$$
M(\phi) \to [0,1]/\left\{0\sim 1\right\} = S^1
$$
whose fiber is $S$ and whose geometric monodromy is $\phi$.

Now suppose that $\phi : (S,x_0) \to (S,x_0)$ is a diffeomorphism of
$S$ that represents an element of $T_g^1$.  The mapping torus bundle
$$
M(\phi) \to S^1
$$
has a canonical section $\sigma : S^1 \to M(\phi)$ which takes
$t \in S^1$ to $(x_0,t)\in M(\phi)$.

Denote the ``jacobian'' of $S$, $H_\dot(S,\R/\Z)$,
by $\Jac S$. The next task is to imbed $M(\phi)$ into $\Jac S$ using
the section $\sigma$ of base points.
To this end, choose a basis $\w_1,\dots, \w_{2g}$ of
$H^1(S,\Z)$. This gives an identification of $\Jac S$ with
$\left(\R/\Z\right)^{2g}$.  Choose closed,
real-valued 1-forms $w_1,\dots,w_{2g}$
that represent $\w_1,\dots, \w_{2g}$. These have integral periods.
Since $\phi$ acts trivially on $H^1(S)$, there are smooth functions
$f_j : S \to \R$ such that
$$
\phi^\ast w_j = w_j + df_j.
$$
These functions are uniquely determined if we insist, as we shall, that
$f_j(x_0)=0$ for each $j$. Set
$$
\wv = (w_1,\dots,w_g)\text{ and } \f = (f_1,\dots, f_g).
$$
The map
$$
S \times [0,1] \to \Jac S
$$
defined by
$$
(x,t) \mapsto t\f(x) + \int_{x_0}^x \wv
$$
preserves the equivalence relations of the mapping torus $M(\phi)$,
and therefore induces a map
$$
\nu(\phi) : (M(\phi),\sigma(S^1)) \to (\Jac S,0).
$$
Define $\taut(\phi)$ to be the homology class of $M(\phi)$ in
$H_3(\Jac S,\Z)$:
$$
\taut(\phi) = \nu(\phi)_\ast[M(\phi)]\in \Lambda^3 H_1(S,\Z).
$$

\begin{proposition} If $\phi$, $\psi$ are diffeomorphisms of $S$ that
act trivially on $H_1(S)$, then
\begin{enumerate}
\item[(a)] $\taut(\phi)$ is independent of the choice of representatives
$w_1,\dots, w_g$ of the basis $\w_1,\dots, \w_{2g}$ of $H^1(S,\Z)$;
\item[(b)] $\taut(\phi)$ is independent of the choice of basis
$\w_1,\dots \w_{2g}$ of $H^1(S,\Z)$;
\item[(c)] $\taut(\phi)$ depends only on the isotopy class of $\phi$;
\item[(d)] $\taut(\phi\psi) = \taut(\phi) + \taut(\psi)$;
\item[(e)] $\taut(g \psi g^{-1}) = g_\ast \taut(\phi)$ for
all diffeomorphisms
$g$ of $S$, where $g_\ast$ is the automorphism of $\Lambda^3 H_1(S)$
induced by $g$.
\end{enumerate}
\end{proposition}

\begin{pf}
If $w'_1,\dots, w'_{2g}$ is another set of representatives of the $\w_j$,
then there are functions $g_j : S \to \R$ such that
$w'_j = w_j + dg_j \text{ and } g_j(x_0) = 0$.
For each $s \in [0,1]$, the 1-form $w_j(s) = w_j + sdg_j$
is closed on $S$ and represents $\w_j$. The map
$$
\nu_s : M(\phi) \to \Jac S
$$
defined using the representatives $w_j(s)$ takes $(x,t)$ to
$$
\left(t\Bigl(f_j(x) + s\bigl(g_j(\phi(x))-g_j(x)\bigr)\Bigr) + sg_j(x) +
\int_{x_0}^x w_j\right).
$$
Since  this depends continuously on $s$, it follows that $\nu_0$ is
homotopic to $\nu_1$. The first assertion follows.

The second assertion follows from linear algebra. The proof of the
third assertion is similar to that of the first.

To prove the fourth assertion, observe that the quotient of
$M(\phi\psi)$ obtained by identifying $(x,1)$ with $(\psi(x),1/2)$ is
the union of $M(\phi)$ and $M(\psi)$. The map $\nu(\phi\psi)$ factors
through the quotient $M(\phi)\cup M(\psi)$ of $M(\phi\psi)$,
and its restrictions to $M(\phi)$ and $M(\psi)$ are
$\nu(\phi)$ and $\nu(\psi)$, respectively. Additivity follows.

Suppose that $g:(S,x_0) \to (S,x_0)$ is a diffeomorphism. The map
$(g,\id) : S\times [0,1] \to S\times [0,1]$ induces a diffeomorphism
$$
F(g) : M(\phi) \to M(g \phi g^{-1}).
$$
To prove the last assertion, it suffices to prove that the
diagram
$$
\begin{CD}
M(\phi)  @>{F(g)}>> M(g \phi g^{-1}) \cr
@V\nu(\phi)VV   @VV\nu(g\phi g^{-1})V \cr
\Jac S @>g_\ast>> \Jac S \cr
\end{CD}
$$
commutes up to homotopy. In the proof of the first assertion, we
saw that the homotopy class of $\nu$ depends only on the basis
of $H^1(S,\Z)$ and not on the choice of de~Rham representatives.
Set $w'_j = g^\ast w_j$. Since the diagram
$$
\begin{CD}
H_1(S) @>g_\ast>> H_1(S) \cr
@V\int w'_jVV @VV\int w_jV\cr
\R @= \R\cr
\end{CD}
$$
commutes, it suffices to prove that the diagram
$$
\begin{CD}
M(\phi)  @>{F(g)}>> M(g \phi g^{-1}) \cr
@V{\nu'}VV   @VV{\nu}V \cr
\left(\R/\Z\right)^{2g} @>\id>>
\left(\R/\Z\right)^{2g} \cr
\end{CD}
$$
commutes, where $\nu$ is defined using $w_1,\dots w_{2g}$, and
$\nu'$ is defined using the representatives $w'_1,\dots w'_{2g}$.
This last assertion is easily verified.
\end{pf}

Recall that the homology groups of $T_g^1$ are $Sp_g(\Z)$ modules;
the action on $H_1(T_g)$ is given by
$$
g : [\phi] \mapsto [\tilde{g}\phi\tilde{g}^{-1}],
$$
where $g\in Sp_g(\Z)$ and $\tilde{g}$ is any element of
$\Gamma_g^1$ that projects to $g$ under the canonical
homomorphism.

\begin{corollary}
The map $\taut$ induces an $Sp_g(\Z)$-equivariant homomorphism
$$
\tau_g^1 : H_1(T_g^1,\Z) \to \Lambda^3 H_1(S,\Z).
$$
\end{corollary}

{}From $\tau_g^1$, we can construct a representation $\tau_g$ of
$H_1(T_g)$. The kernel of the natural surjection $T_g^1 \to T_g$
is isomorphic to $\pi_1(S,x_0)$. The composition of the induced map
$H_1(S,\Z) \to H_1(T_g^1,\Z)$ with $\tau_g^1$ is  easily seen to be the
canonical inclusion
$$
\blank\times [S] : H_1(S,\Z) \hookrightarrow H_3(\Jac S,\Z)
$$
induced by taking Pontrjagin product with $\nu_\ast[S]$. We therefore
have an induced $Sp_g(\Z)$-equivariant map
$$
\tau_g : H_1(T_g,\Z) \to \Lambda^3 H_1(S,\Z) /H_1(S,\Z).
$$

The following result of Johnson is an immediate corollary of
Theorem \ref{john:h1}.

\begin{theorem}\label{tau_g}
The homomorphism $\tau_g$ is an isomorphism modulo 2-torsion.
\end{theorem}

It is not difficult to boot strap up from Johnson's basic computation
to prove the following result.

\begin{theorem}\label{h1_alg}
There is a natural $Sp_g(\Z)$-equivariant isomorphism
$$
\tau_{g,r}^n : H_1(T_{g,r}^n,\Q) \to H_1(S,\Q)^{\oplus(n+r)} \oplus
\Lambda^3 H_1(S,\Q) /H_1(S,\Q).
$$
\end{theorem}

An important consequence of Johnson's theorem is that the action
of $Sp_g(\Z)$ on $H_1(T_{g,r}^n,\Q)$ factors through a rational
representation of the $\Q$-algebraic group $Sp_g$. Let $\lambda_1,
\dots, \lambda_g$ be a fundamental set of dominant integral weights
of $Sp_g$. Denote the irreducible $Sp_g$-module with highest weight
$\lambda$ by $V(\lambda)$. The fundamental representation of $Sp_g$ is
$H_1(S)$. It is well known (and easily verified) that
$$
\Lambda^3 H_1(S) \approx V(\lambda_1) \oplus V(\lambda_3).
$$
The previous result can be restated by saying that
$$
H_1(T_{g,r}^n,\Q) \approx
V(\lambda_3) \oplus V(\lambda_1)^{\oplus(n+r)}
$$
as $Sp_g$ modules.

\section{A Second Definition of the Johnson Homomorphism}
\label{original}

In this section we relate the definition of $\tau_g^1$ given in the
previous section to Johnson's original definition, which is defined
using the action of $T_g^1$ on the lower central series of
$\pi_1(S,x_0)$. It is better suited to computations.  In order to
relate this definition to the one given in the previous section,
we need to study the cohomology ring of the mapping torus associated
to an element of the Torelli group.

Suppose that the diffeomorphism $\phi : (S,x_0) \to (S,x_0)$
represents an element of $T_g^1$. As explained in the previous
section, the associated mapping torus $M=M(\phi)$ fibers over $S^1$
and has a canonical section $\sigma$. This data guarantees that there
is a canonical decomposition of the cohomology of $M$.

Since $\phi$ acts trivially on the homology of $S$, the $E_2$-term of
the Leray-Serre spectral sequence of the fibration $\pi : M \to S^1$
satisfies
$$
E_2^{r,s} = H^r(S^1)\otimes H^s(S).
$$
This spectral sequence degenerates for trivial reasons. Consequently,
there is a short exact sequence
$$
0 \to H^1(S^1,\Z) \stackrel{\pi^\ast}{\to} H^1(M,\Z) \stackrel{i^\ast}
{\to} H^1(S,\Z) \to 0,
$$
where $\pi$ is the projection to $S^1$ and $i : S \hookrightarrow M$
is the inclusion of the fiber over the base point $t=0$ of $S^1$.
The section $\sigma$ induces a splitting of this sequence. Denote
$\pi^\ast$ of the positive generator of $H^1(S^1,\Z)$ by $\theta$. Then
we have the decomposition
\begin{equation}\label{h1}
H^1(M,\Z) = H^1(S,\Z) \oplus \Z\theta.
\end{equation}
{}From the spectral sequence, it follows that we have an exact sequence
$$
0 \to \theta\wedge H^1(S,\Z) \to H^2(M,\Z)
\stackrel{i^\ast}{\to} H^2(S,\Z) \to 0.
$$
Denote the Poincar\'e dual of a homology class $u$ in $M$ by $PD(u)$.
Since
$$
\int_S PD(\sigma) = \sigma\cdot S = 1
$$
it follows that the previous sequence can be split by taking the
positive generator of $H^2(S,\Z)$ to $PD(\sigma)$. We therefore
have a canonical splitting
\begin{equation}\label{h2}
H^2(M,\Z) = \Z PD(\sigma) \oplus \theta\wedge H^1(S,\Z).
\end{equation}

The cup product pairing
$$
c : H^1(M)\otimes H^2(M) \to H^3(M) \approx \Z
$$
induces pairings between the summands of the decompositions
(\ref{h1}) and (\ref{h2}).

\begin{proposition}\label{prop}
The cup product $c$ satisfies:
\begin{enumerate}
\item[(a)] $c(\theta \otimes PD(\sigma)) = 1$;
\item[(b)] the restriction of $c$ to $H^1(S) \otimes PD(\sigma)$
vanishes;
\item[(c)] the restriction of $c$ to
$\theta \otimes \bigl(\theta\wedge H^1(S)\bigr)$ vanishes;
\item[(d)] the restriction of $c$ to
$H^1(S) \otimes \bigl( \theta\wedge H^1(S)\bigr)$ takes
$u \otimes (\theta\wedge v)$ to $-\int_S u\wedge v$.
\end{enumerate}
\end{proposition}

\begin{pf}
Since $\theta$ is the Poincar\'e dual of the fiber $S$, we have
$$
\int_M \theta \wedge PD(\sigma) = \int_M PD(S) \wedge PD(\sigma)
= S\cdot \sigma = 1.
$$

In the decomposition (\ref{h1}), $H^1(S)$ is identified with the kernel
of $\sigma^\ast: H^1(M) \to H^1(S^1)$; that is, with those
$u\in H^1(M)$ such that
$$
\int_\sigma u = 0.
$$
The second assertion now follows as
$$
\int_M u \wedge PD(\sigma) = \int_\sigma u
$$
for all $u \in H^1(M)$. The third and fourth assertions are easily
verified.
\end{pf}

To complete our understanding of the cohomology ring of $M$, we
consider the cup product
$$
\Lambda^2 H^1(M) \to H^2(M).
$$
Since $\theta\wedge \theta =0$, there is only one interesting part of
this mapping; namely, the component
$$
\Lambda^2 H^1(S) \to \Z PD(\sigma) \oplus \theta\wedge H^1(S).
$$
There is a unique function
$$
f_\phi : \Lambda^2 H^1(S,\Z) \to H^1(S,\Z)
$$
such that
$$
u\wedge v \mapsto
\left(\int_S u\wedge v,-\theta\wedge f_\phi(u\wedge v)\right)
\in H^2(M,\Z)
$$
with respect to the decomposition (\ref{h2}).
We can view $f_\phi$ as an element of
$$
H_1(S,\Z) \otimes \Lambda^2 H^1(S,\Z).
$$
Using Poincar\'e duality on the last two factors, $f_\phi$ can be
regarded as an element $F(\phi)$ of
$$
H_1(S,\Z)\otimes \Lambda^2 H_1(S,\Z).
$$

There is a canonical imbedding of $\Lambda^3 H_1(S,\Z)$ into
this group. It is defined by
$$
a\wedge b \wedge c \mapsto a\otimes (b\wedge c) +
b\otimes (c\wedge a) + c\otimes (a\wedge b).
$$

\begin{theorem}\label{image}
The invariant $F(\phi)$ of the cohomology ring of $M(\phi)$ is the
image of $\taut(\phi)$ under the canonical imbedding
$$
\Lambda^3 H_1(S,\Z) \hookrightarrow
H_1(S,\Z)\otimes \Lambda^2 H_1(S,\Z).
$$
\end{theorem}

\begin{pf}
The dual of $\tau_g^1(\phi)$ is the map
$$
\Lambda^3 H^1(S) \to \Z
$$
defined by
$$
u \wedge v \wedge w \mapsto \int_{M(\phi)} u \wedge v \wedge w.
$$
Here we have identified $H^1(S)$ with $H^1(\Jac S)$ using the
canonical isomorphism
$$
\nu^\ast : H^1(\Jac S) \stackrel{\sim}{\to} H^1(S).
$$
The map $\nu(\phi) : M \to \Jac S$ collapses $\sigma$ to the point 0.
It follows that the image
of
$$
\nu(\phi)^\ast : H^1(\Jac S) \to H^1(M)
$$
lies in the subspace we are identifying with $H^1(S)$ in the
decomposition (\ref{h1}). Since the restriction of $\nu(\phi)$ to
the fiber over the base point $t=0$ of $S^1$ is the isomorphism
$\nu^\ast$, it follows that the diagram
$$
\begin{CD}
H^1(\Jac S) @>{\nu(\phi)^\ast}>> H^1(M) \cr
@V{\nu^\ast}VV  @| \cr
H^1(S) @>{i}>>  H^1(M) \cr
\end{CD}
$$
commutes, where $i$ is the inclusion given by the splitting
(\ref{h1}). That is, all the identifications we have made with $H^1(S)$
are compatible.

We will compute the dual of $\tau_g^1(\phi)$ using $F(\phi)$,
which we regard as a homomorphism
$$
F(\phi) : H^1(S) \otimes \Lambda^2 H^1(S) \to \Z
$$
It follows from (\ref{prop}) that this map takes
$u\otimes (v\wedge w)$ to
$$
\int_S u\wedge f_\phi(v\wedge w).
$$
The assertion that $F(\phi)$ lie in $\Lambda^3H_1(S)$ is equivalent
to the assertion that
$$
F(\phi)(u\otimes (v\wedge w)) = F(\phi)(v\otimes (w\wedge u))
=F(\phi)(w\otimes (u\wedge v)),
$$
which is easily verified using (\ref{prop}). The equality of
$F(\phi)$ and $\tau_g^1(\phi)$ follows as
$$
\tau_g^1(\phi)(u\wedge v\wedge w)
= \int_{M} u\wedge v\wedge w
= -\int_{M} u \wedge \theta \wedge f_\phi(v\wedge w)
=  F(\phi) (u\otimes (v\wedge w)).
$$
\end{pf}

We are now ready to give Johnson's original definition of $\tau_g^1$.
Denote the lower central series of a group $\pi$ by
$$
\pi = \pi^{(1)} \supseteq \pi^{(2)} \supseteq \pi^{(3)}
\supseteq \cdots
$$
We regard the cup product
$$
\Lambda^2 H^1(S,\Z) \to H^2(S,\Z) \approx \Z
$$
as an element $q$ of $\Lambda^2 H^1(S,\Z)$.

\begin{proposition}\label{lcs}
The commutator mapping
$$
\com : \pi_1(S,x_0) \times \pi_1(S,x_0) \to \pi_1(S,x_0)
$$
induces an isomorphism
$\Lambda^2 H_1(S,\Z)/q \to \pi_1(S,x_0)^{(2)} /
\pi_1(S,x_0)^{(3)}$.
\end{proposition}

\begin{pf}
This follows directly from the standard fact (see \cite{serre}
or \cite{mks}) that if $F$ is a free group, the commutator
induces an isomorphism
$$
\Lambda^2 H_1(F) \stackrel{\sim}{\to} F^{(2)}/F^{(3)}
$$
and from the standard presentation of $\pi_1(S,x_0)$.
\end{pf}

An element $\phi$ of $T_g^1$ induces an automorphism of
$\pi_1(S,x_0)$. Since it acts trivially on $H_1(S)$,
$$
\phi(\gamma)\gamma^{-1} \in \pi_1(S,x_0)^{(2)}
$$
for all $\gamma\in \pi_1(S,x_0)$. From (\ref{lcs}), it follows
that $\phi$ induces a well defined map
$$
\tauh(\phi) : H_1(S,\Z) \to \Lambda^2 H_1(S,\Z)/q
$$
Using Poincar\'e duality, we may view this as an element
$L(\phi)$ of
$$
H_1(S,\Z) \otimes \left(\Lambda^2 H_1(S,\Z)/q\right).
$$

\begin{theorem}
The image of $F(\phi)$ in
$H_1(S,\Z) \otimes \left(\Lambda^2 H_1(S,\Z)/q\right)$
is $L(\phi)$.
\end{theorem}

\begin{pf}
Since $H_1(S,\Z) \otimes \left(\Lambda^2 H_1(S,\Z)/q\right)$
is torsion free, it suffices to show that the image of $F(\phi)$
in
$$
H_1(S,\Q) \otimes \left(\Lambda^2 H_1(S,\Q)/q\right)
$$
is $L(\phi)$. For the rest of this proof, all (co)homology groups have
$\Q$ coefficients.

For all groups $\pi$ with finite dimensional $H_1(\phantom{x},\Q)$, the
sequence
\begin{equation}\label{ex-seq}
0 @>>> H^1(\pi)  @>h^\ast>>
\left(\pi/\pi^{(3)}\right)^\ast @>{\com^\ast}>>
\Lambda^2 H^1(\pi) @>{\wedge}>> H^2(\pi)
\end{equation}
of $\Q$ vector spaces is exact. Here $(\phantom{x})^\ast$ denotes
the dual vector space, $h^\ast$ the dual Hurewicz homomorphism, and
$\com$ the map induced by the commutator. This can be proved
using results in either \cite[\S 2.1]{chen} or
\cite[\S 8]{sullivan:inf}.

We apply this sequence to the fundamental group
of the mapping torus. Choose $m_0=(x_0,0)$ as the base point of
$M$.  Since $M$ fibers over the circle with fiber $S$,
we have an exact sequence
$$
1 \to \pi_1(S,x_0) \to \pi_1(M,m_0) \to \Z \to 0.
$$
The section $\sigma$ induces a splitting $\Z \to \pi_1(M,m_0)$.

Denote the image of $1$ by $\sigma$. Observe that if
$\gamma\in \pi_1(S,x_0)$, then
$$
\sigma \gamma \sigma^{-1} = \phi(\gamma).
$$
It follows that the inclusion $\pi_1(S,x_0) \hookrightarrow
\pi_1(M,m_0)$
induces isomorphisms
$$
\pi_1(S,x_0)^{(k)} \approx \pi_1(M,m_0)^{(k)}
$$
for all $k> 1$ and, as above, that $\sigma$ induces an
isomorphism
$$
H_1(M) = H_1(S) \oplus \Q\Sigma,
$$
where $\Sigma$ denotes the homology class of $\sigma$.
It also follows that for all $a\in H_1(S)$
$$
\taut(\phi)(a) = [\Sigma,a] \in \pi_1(M)^{(2)}/\pi_1(M)^{(3)}
\approx \pi_1(S)^{(2)}/\pi_1(S)^{(3)}.
$$

Using (\ref{prop}) and the exact sequence (\ref{ex-seq}), we
see that for all $u\in H_1(S)$, the image of $f_\phi^\ast (u)$ in
$$
\Lambda^2 H_1(S)/q
$$
is $[\Sigma,u]$, which is $\tauh(\phi)(u)$ as we have seen. The
result follows.
\end{pf}

The composite of the inclusion
$$
\Lambda^3 H_1(S,\Z) \hookrightarrow H_1(S,\Z)\otimes \Lambda^2 H_1(S,\Z)
$$
with the quotient mapping
$$
H_1(S,\Z)\otimes \Lambda^2 H_1(S,\Z) \to
H_1(S,\Z)\otimes \left(\Lambda^2 H_1(S,\Z)/q\right)
$$
is injective. One way to see this is to tensor with $\Q$ and
note that both of these maps are maps of $Sp_g$ modules. One can
then use the fact that $\Lambda^3 H_1(S)$ is the sum of the first
and third fundamental representations of $Sp_g$ to check the result.
The following result is therefore a restatement of (\ref{image}).

\begin{corollary}\label{equality}
$L(\phi)$ lies in the image of the canonical injection
$$
\Lambda^3 H_1(S,\Z) \hookrightarrow
H_1(S,\Z)\otimes \left(\Lambda^2 H_1(S,\Z)/q\right)
$$
and the corresponding point of $\Lambda^3 H_1(S)$ is
$\tau_g^1(\phi)$.
\end{corollary}

In his fundamental papers, Johnson defines $\tau_g^1(\phi)$ to
be $L(\phi)$. The other two definitions we have given were stated
in \cite{johnson:survey}.

\section{Picard Groups}

In \cite{mumford}, Mumford showed that
$$
c_1 : \Pic \M_g \otimes \Q\to H^2(\M_g,\Q)
$$
is an isomorphism. Using Johnson's computation of $H_1(T_g,\Q)$ and
the well known result (\ref{pic}), we will prove the analogous
statement for all $\M_{g,r}^n(L)$ when $g\ge 3$. The novelty
lies in the variation of the level, and not in the variation
of the decorations $r$ and $n$. The first, and principal, step is to
establish the vanishing of the $H^1(\M_{g,r}^n(L))$.

\begin{proposition}\label{h1-level}
Suppose that $L$ is a finite index subgroup of $Sp_g(\Z)$.
If $g \ge 3$, then $H^1(\M_{g,r}^n(L),\Z)=0$.
\end{proposition}

Since $H^1(\phantom{X},\Z)$ is always torsion free, it suffices to
prove that $H^1(\M_g(L),\Q)$ vanishes.  We will prove a stronger
result.

\begin{proposition}\label{calc}
Suppose that $L$ is a finite index subgroup of $Sp_g(\Z)$ and that
$g\ge 3$. If $V(\lambda)$ is an irreducible representation of $Sp_g$
with highest weight $\lambda$, then
$$
H^1(\Gamma_{g,r}^n(L),V(\lambda)) =
\begin{cases}
\Q^{r+n} & \text{if } \lambda = \lambda_1; \cr
\Q & \text{if }\lambda = \lambda_3;\cr
0 & otherwise.
\end{cases}
$$
Consequently, $H^1(\M_{g,r}^n(L),\Z)$ vanishes for all $r,n$ when
$g\ge 3$.
\end{proposition}

\begin{pf}
It follows from the Hochschild-Serre spectral sequence
$$
H^r(L,H^s(T_{g,r}^n\otimes V(\lambda))) \implies
H^{r+s}(\Gamma_{g,r}^n(L),V(\lambda))
$$
that there is an exact sequence
\begin{multline*}
0 \to H^1(L,V(\lambda)) \to
H^1(\Gamma_{g,r}^n(L),V(\lambda)) \\
\to H^0(L,H^1(T_{g,r}^n\otimes V(\lambda)))
@>{d_2}>> H^2(L,V(\lambda)).
\end{multline*}
By a result of Ragunathan
\cite{ragunathan}, the first term vanishes when $g\ge 2$. By
(\ref{h1_alg}), the third term vanishes except when $\lambda$ is
either $\lambda_1$ or $\lambda_3$. This proves the result except
in the cases when $\lambda$ is either $\lambda_1$ or $\lambda_3$.
In these exceptional cases, the third term has
rank $r+n$ or 1, respectively. To complete the proof, we need to show
that the differential $d_2$ is zero.

There are several ways to do this. Perhaps the most straightforward is
to use the result, due to Borel \cite{borel:twisted}, which asserts
that the last group vanishes when $g\ge 8$. This establishes the result
when $g\ge 8$. When $r\ge 1$, the vanishing of $d_r$ for all $g\ge 3$
follows as the diagram
$$
\begin{CD}
H^0(L,H^1(T_{g,r}^n\otimes V(\lambda))) @>d_2>> H^2(L,V(\lambda)) \cr
@AAA						@AAA \cr
H^0(L,H^1(T_{g+8,r}^n\otimes V(\lambda))) @>d_2>> H^2(L_{g+8},V(\lambda))\cr
\end{CD}
$$
commutes. Here $L_{g+8}$ is any finite index subgroup of $Sp_{g+8}(\Z)$
such that
$$
L_{g+8}\cap Sp_g(\Z) \subseteq L
$$
and the vertical maps are induced by the ``stabilization map''
$$
\Gamma_{g,r}^n(L) \to \Gamma_{g+8,r}^n(L_{g+8}).
$$
When $r=0$ and $\lambda=\lambda_1$, there is nothing to prove. This
leaves only the case $r=0$ and $\lambda = \lambda_3$ which follows
as the diagram
$$
\begin{CD}
H^0(L,H^1(T_g^n\otimes V(\lambda))) @>d_2>> H^2(L,V(\lambda)) \cr
@AAA						@AAA \cr
H^0(L,H^1(T_{g,1}^n\otimes V(\lambda))) @>d_2>> H^2(L,V(\lambda))\cr
\end{CD}
$$
which arises from the homomorphism $\Gamma_{g,1}^n \to \Gamma_g^n$,
commutes.
\end{pf}

Denote the category of $\Z$ mixed Hodge structures by $\H$.
We shall denote the group of ``integral $(0,0)$ elements''
$\Hom_\H(\Z,H)$ of a mixed Hodge structure $H$ by $\Gamma H$.

Suppose that $X$ is a smooth variety. Since $H^1(X,\Z)$ is torsion
free, we can define
$$
W_1H^1(X,\Z) = W_1H^1(X,\Q) \cap H^1(X,\Z).
$$
This is a polarized, torsion free Hodge structure of weight 1. Set
$$
JH^1(X) = {W_1H^1(X,\C) \over W_1H^1(X,\Z) + F^1W_1H^1(X,\C)}.
$$
This is a polarized Abelian variety.

\begin{theorem}\label{pic}
If $X$ is a smooth variety, then there is a natural exact sequence
$$
0 \to JH^1(X) \to \Pic X \to \Gamma H^2(X,\Z(1)) \to 0.
$$
\end{theorem}

Alternatively, this theorem may be stated as saying that the
cycle map
$$
\Pic X \to H_\H^2(X,\Z(1))
$$
is an isomorphism, where $H_\H^\dot$ denotes Beilinson's absolute
Hodge cohomology, the refined version of Deligne cohomology
defined in \cite{beilinson:ahc}.

\begin{pf}
Choose a smooth completion $\X$ of $X$ for which $\X-X$
is a normal crossings divisor $D$ in $\X$ with smooth components.
Denote the dimension of $X$ by $d$. From the usual exponential
sequence, we have a short exact sequence
$$
0 \to JH^1(\X) \to \Pic \X \to \Gamma H^2(X,\Z) \to 0.
$$
{}From \cite[(1.8)]{fulton}, we have an exact sequence
$$
CH^0(D) \to \Pic \X \to \Pic X \to 0.
$$
The Gysin sequence
$$
0 \to H^1(\X) \to H^1(X) \to H_{2d-2}(D)(-2d) \to H^2(\X) \to H^2(X)
\to H_{2d-3}(D)(-2d)
$$
is an exact sequence of $\Z$ Hodge structures. Since $H_{2d-2}(D)(-2d)$ is
torsion free and of weight 2, it follows that
$$
W_1H^1(X,\Z) = H^1(\X,\Z),
$$
and therefore that $JH^1(X) = JH^1(\X)$.
Next, since each component $D_i$ of $D$ is smooth, it follows that
$$
H_{2d-3}(D)(-2d) = \bigoplus_i H^1(D_i,\Z)(-1),
$$
and is therefore torsion free and of weight 3. It follows that the
sequence
$$
H_{2d-2}(D)(-2d) \to \Gamma H^2(\X) \to \Gamma H^2(X) \to 0
$$
is exact. Since the cycle map
$$
CH^0(D) \to H_{2d-2}(D)
$$
is an isomorphism \cite[(1.5)]{fulton}, the result follows.
\end{pf}

It is now an easy matter to show that the Picard groups of the
$\M_{g,r}^n(L)$ are finitely generated.

\begin{theorem}
Suppose that $L$ is a finite index subgroup of $Sp_g(\Z)$.
If $g\ge 3$, then for all $r,n$, the Chern class map
$$
c_1 : \Pic \M_{g,r}^n(L) \to \Gamma H^2(\M_{g,r}^n(L),\Z)
$$
is an isomorphism when $\Gamma_{g,r}^n(L)$ is torsion free, and is
an isomorphism after tensoring with $\Q$ in general.
\end{theorem}

\begin{pf}
The case when $\Gamma_{g,r}^n(L)$ is torsion free follows directly
from (\ref{h1-level}) and (\ref{pic}). To prove the assertion in
general, choose a finite index normal subgroup $L'$ of $L$ such that
 $\Gamma_{g,r}^n(L')$ is torsion free. Let
$$
G = \Gamma_{g,r}^n(L) / \Gamma_{g,r}^n(L') \approx L/L'.
$$
Then it follows from the Teichm\"uller description of moduli spaces
that the projection
$$
\pi : \M_{g,r}^n(L') \to \M_{g,r}^n(L)
$$
is a Galois covering with Galois group $G$. It follows from the first
case that
$$
c_1 :\Pic \M_{g,r}^n(L') \to \Gamma H^2(\M_{g,r}^n(L'),\Z)
$$
is a $G$-equivariant isomorphism. The result now follows as the
projection $\pi$ induces isomorphisms
$$
\Pic \M_{g,r}^n(L)\otimes\Q \approx
H^0(G,\Pic \M_{g,r}^n(L')\otimes \Q)
$$
and
$$
\Gamma H^2(\M_{g,r}^n(L),\Q) \approx \Gamma H^0(G,H^2(\M_{g,r}^n(L'),\Q)).
$$
\end{pf}

If we knew that $H^2(T_g,\Q)$ were finite dimensional and a rational
representation of $Sp_g$, we would know from Borel's work
\cite{borel:triv} that $H^2(\M_{g,r}^n(L),\Q)$ would be independent
of the level $L$, once $g$ is sufficiently large; $g\ge 8 $ should
do it --- cf.\ \cite{borel:twisted}. It would then follow, for
sufficiently large $g$, that the Picard number
of $\M_{g,r}^n(L)$ is $n+r+1$. At present it is not even known whether
$H^2(T_g,\Q)$ is finite dimensional. The computation
of this group, and the related problem of finding a presentation of
$T_g$ appear to be deep and difficult. It should be mentioned that
the only evidence for the belief that the Picard number of each
$\M_g(L)$ is one comes from Harer's computation \cite{harer:spin}
of the Picard numbers of the moduli spaces of curves with a
distinguished theta characteristic.

\section{Normal Functions}

In this section, we define abstract normal functions which generalize
the normal functions of Poincar\'e and Griffiths.
We begin by reviewing how a family of homologically trivial
algebraic cycles in a family of smooth projective varieties gives
rise to a normal function.

Suppose that $X$ is a smooth variety. A homologically trivial
algebraic $d$-cycle in $X$ canonically determines an element of
$$
\Ext^1_\H(\Z,H_{2d+1}(X,\Z(-d))).
$$
This extension is obtained by pulling back the exact sequence
$$
0 \to H_{2d+1}(X,\Z(-d)) \to H_{2d+1}(X,Z,\Z(-d)) \to H_{2d}(Z,\Z(-d))
\to \cdots
$$
of mixed Hodge structures along the inclusion
$$
\Z \to H_{2d}(|Z|,\Z(-d))
$$
that takes 1 to the class of $Z$.

When $H$ is a mixed Hodge structure all of whose weights are
non-positive, there is a natural isomorphism
$$
J H \approx \Ext^1_\H(\Z,H),
$$
where
$$
JH = {H_\C \over F^0H_\C + H_\Z}.
$$
(This is well known---see for example \cite{carlson}. Our
conventions will be taken from \cite[(2.2)]{hain:heights}.)

When $X$ is projective,
Poincar\'e duality provides an isomorphism of the complex torus
$JH_{2d+1}(X,\Z(-d))$ with the Griffiths intermediate jacobian
$$
\Hom_\C(F^dH^{d+1}(X),\C)/H_{2d+1}(X,\Z).
$$
The point in $JH_{2d+1}(X,Z(d))$ corresponding to the cycle $Z$
under this isomorphism is $\int_\Gamma$, where $\Gamma$ is a
real $2d+1$ chain that satisfies $\partial \Gamma = Z$.

Now suppose that $\calX \to T$ is a family of smooth projective varieties
over a smooth base $T$. Suppose that $\calZ$ is an algebraic cycle in
$\calX$ which is proper over $T$ of relative dimension $d$. Denote the
fibers of $\calX$ and $\calZ$ over $t\in T$ by $X_t$ and $Z_t$, respectively.

The set of $H_{2d+1}(X_t,\Z(-d))$ form a variation of Hodge structure
$\V$ over $T$ of weight $-1$. We can form the relative intermediate
jacobian
$$
\J_d \to T;
$$
this has fiber $JH_{2d+1}(X_t,\Z(-d))$ over $t\in T$. The family of
cycles $\calZ$ defines a section this bundle. Such a section is what
Griffiths calls {\it the normal function of the cycle $\calZ$}
\cite{griffiths}. Griffiths' normal functions generalize those of
Poincar\'e.

We will generalize this notion further. Before we do, note that the
elements of $\Ext^1_\H(\Z,H_{2d+1}(X_t,\Z(-d)))$ defined by the
cycles $Z_t$ fit
together to form a variation of mixed Hodge structure over $T$. It
follows from the main result of \cite{g-n-p} that this variation  is
good in the sense of \cite{steenbrink-zucker} along each curve in $T$,
and is therefore good in the sense of Saito \cite{saito}.

Suppose that $T$ is a smooth variety and that $\V \to T$ is a
variation of Hodge structure over $T$ of negative weight.  Denote the
bundle over $T$ whose fiber over $t\in T$ is
$$
JV_t \approx \Ext^1_\H(\Z,V_t)
$$
by $J\calV$.

\begin{definition}
A holomorphic section $s: T \to J\calV$ of $J\calV \to T$ is a {\it normal
function} if it defines an extension
$$
0 \to \V \to \E \to \Z_T \to 0
$$
in the category $\H(T)$ of good variations of mixed Hodge structure
over $T$.
\end{definition}

\begin{remark}\label{norm}
We know from the preceding discussion that families of
homologically trivial cycles in a family $\calX \to T$ define normal
functions in this sense.
\end{remark}

The asymptotic properties of good variations of mixed Hodge structure
guarantee that these normal functions have nice properties.

\begin{lemma}[Rigidity]\label{rigid}
If $\V \to T$ and $\V' \to T$
are two good variations of mixed Hodge structure over $T$ with the
same fiber $V_{t_0}$ (viewed as a mixed Hodge structure) over some
point $t_0$ of $T$ and with the same monodromy representations
$$
\pi_1(T,t_0) \to \Aut V_{t_0},
$$
then $\V_1$ and $\V_2$ are isomorphic as variations.
\end{lemma}

\begin{pf}
The proof is a standard application of the theorem of the fixed part.
The local system $\Hom_\Z(\V,\V')$
underlies a good variation of mixed Hodge structure. From Saito's
work \cite{saito}, we know that each cohomology group of a variety
with coefficients in a good variation of mixed Hodge structure has a
natural mixed Hodge structure. So, in particular,
$$
H^0(T,\Hom_\Z(\V,\V'))
$$
has a mixed Hodge structure, and the restriction map
$$
H^0(T,\Hom_\Z(\V,\V')) \to \Hom_\Z(V_{t_0},V_{t_0}')
$$
is a morphism. The result now follows as there are natural
isomorphisms
$$
H^0(T,\Hom_\Z(\V,\V')) \approx
\Hom_{\Z\pi_1(T,t_0)}(V_{t_0},V_{t_0}')
$$
and
$$
\Gamma H^0(T,\Hom_\Z(\V,\V')) \approx \Hom_{\H(T)}(\V,\V'),
$$
where $\H(T)$ denotes the category of good variations of mixed
Hodge structure over $T$.
\end{pf}

\begin{corollary}\label{rigidity}
Two normal functions $s_1, s_2 : T \to J\calV$ are equal if and only
if there is a point $t_0 \in T$ such that $s_1(t_0) = s_2(t_0)$ and
such that the two induced homomorphisms
$$
(s_j)_\ast : \pi_1(T,t_0) \to \pi_1(J\calV,s_1(t_0))
$$
are equal. \qed
\end{corollary}

\section{Extending Normal Functions}
\label{extending}

The strong asymptotic properties of variations of mixed Hodge
structure imply that almost all normal functions extend across
subvarieties where the original variation of Hodge structure is
non-singular. Suppose that $X$ is a smooth variety and that
$\V$ is a variation of Hodge structure over $X$ of negative
weight. Denote the associated intermediate jacobian bundle by
$\J \to X$.

\begin{theorem}\label{extend}
Suppose that $U$ is a Zariski open subset of $X$ and that
$s: U \to \J|_U$ is a normal function defined on $U$. If
the weight of $\V$ is not $-2$, then $s$ extends to a normal
function $\stilde : X \to \J$.
\end{theorem}

\begin{pf}
Write $U=X-Z$. By Hartog's Theorem, it suffices to show that $s$
extends to a normal function on the complement of the union of the
singular locus of $Z$ and the union of the components of $Z$ of
codimension $\ge 2$ in $X$. That is,  we may assume that $Z$ is a
smooth divisor.

The problem of extending $s$ is local. By taking
a transverse slice, we can reduce to the case where $X$ is the
unit disk $\Delta$ and $Z$ is the origin.
In this case, we have a variation of Hodge structure over $\Delta$.
The normal function $s: \Delta^\ast \to \J$ corresponds to a good
variation of mixed Hodge structure $\E$ over the punctured disk
$\Delta^\ast$ which is an extension
$$
0 \to \V|_{\Delta^\ast} \to \E \to \Z_{\Delta^\ast}.
$$
To prove that the normal function extends, it suffices to show
that the monodromy of $\E$ is trivial, for then the local system
$\E$ extends uniquely as a flat bundle to $\Delta$ and the Hodge
filtration extends across the origin as $\E$ is a good
variation.

Since $\V$ is defined on the whole disk, it has trivial monodromy.
It follows that the local monodromy operator $T$ of $\E$ satisfies
$$
(T-I)^2 = 0
$$
and that the local monodromy logarithm $N$ is $T-I$. Since
$E$ is a good variation, it has a relative weight filtration $M_\dot$
(cf.\  \cite{steenbrink-zucker}), which is defined over $\Q$ and
satisfies $NM_l \subseteq M_{l-2}$. From the defining properties of
$M_\dot$ (\cite[(2.5)]{steenbrink-zucker}), we have
$$
M_0 = \E,\, M_{m} = \V, \text{ and } M_{m-1} = 0,
$$
where $m$ is the weight of $\V$.

In the case where $m=-1$, the proof that $N=0$ is simpler. Since this
case is the most important (as it is the one that applies to normal
functions of cycles), we prove it first. The condition $m=-1$ implies
that $M_{-2}=0$. Since $NM_0 \subseteq M_{-2}$, it follows that $N=0$
and consequently, that the normal function extends.

In general, we use defining property
\cite[(3.13)iii]{steenbrink-zucker} of good variations of mixed Hodge
structure which says that
$$
(E_t, M_\dot, F^\dot_{\lim})
$$
is a mixed Hodge structure and $N$ is an morphism of mixed Hodge
structures of type $(-1,-1)$, where $F^\dot_{\lim}$ denotes the
limit Hodge filtration. In this case, $N$ induces a morphism
$$
\Z \approx Gr^M_0 \to Gr^M_{-2},
$$
which is zero if $m\neq -2$. Since $N$ is a morphism of mixed Hodge
structures, the vanishing of this map implies the vanishing of $N$.
\end{pf}

When $m=-2$, there are normal functions that don't extend. For example,
if we take $\V = \Z(1)$, then the bundle of intermediate jacobians is
the bundle $X\times \C^\ast$ and the normal functions are precisely
the invertible regular functions $f: X \to \C$ --- for details see,
for example, \cite[(9.3)]{hain:poly}.

\section{Normal Functions over $\M_{g,r}^n(L)$}
\label{normal_functions}

Throughout this section, we will assume that $g\ge 3$ and $L$ is a finite
index subgroup of $Sp_g(\Z)$ such that $\Gamma_{g,r}^n(L)$ is torsion free.
With this condition on $L$, $\M_{g,r}^n(L)$ is smooth. Each irreducible
representation of $Sp_g$ defines a polarized  $\Q$ variation of Hodge
structure over $\M_{g,r}^n(L)$ which is unique up to Tate twist --- cf.\
(\ref{exist-unique}). It follows that every rational representation of
$Sp_g$ underlies a polarized $\Z$ variation of Hodge structure over
$\M_g(L)$.

\begin{lemma}\label{fin_gen}
If $\V \to \M_{g,r}^n(L)$ is a good variation of Hodge
structure of negative weight whose monodromy representation
$$
\Gamma_{g,r}^n(L) \to \Aut V_\circ\otimes \Q
$$
factors through a rational representation of $Sp_g$ and contains no
copies of the trivial representation, then the group of
normal functions $s : \M_{g,r}^n(L) \to J\calV$ is finitely generated
of rank bounded by
$$
\dim H^1(\Gamma_{g,r}^n(L),V_\Z).
$$
\end{lemma}

\begin{pf}
A normal function corresponds to a variation of mixed Hodge
structure whose underlying local system is an extension
$$
0 \to \V \to \E \to \Z \to 0
$$
of the trivial local system by $\V$.

One can form the semidirect product $\Gamma_{g,r}^n(L)\ltimes V_\Z$,
where the mapping class group acts on $V_\Z$ via a representation
$L\to \Aut V$. The monodromy representation of the local system
$\E$ gives a splitting
$$
\rho : \Gamma_{g,r}^n(L) \to \Gamma_{g,r}^n(L)\ltimes V_\Z
$$
of the natural projection
\begin{equation}\label{mono}
\Gamma_{g,r}^n(L)\ltimes V_\Z \to \Gamma_{g,r}^n(L).
\end{equation}
The splitting is well defined up to conjugation by an element of $V_\Z$.

The first step in the proof is to show that an extension of $\Q$ by $\V$
in the category of $\Q$ variations of mixed Hodge structure is determined
by its monodromy representation. Two such variations can
be regarded as elements of the group
\begin{equation}\label{ext}
\Ext^1_{\H(\M_{g,r}^n(L))}(\Q,\V).
\end{equation}
It is easily seen that their difference is an extension whose monodromy
representation factors through the homomorphism
$\Gamma_{g,r}^n \to Sp_g(\Q)$. It now follow from (\ref{decomp}) and the
assumption that $\V$ contain no copies of the trivial representation
that this difference is the trivial element of (\ref{ext}). The
assertion follows.

{}From \cite[p.~106]{maclane} it follows that the set of splittings of
(\ref{mono}), modulo conjugation by elements of $V_\Z$, is isomorphic
to
$$
H^1(\Gamma_{g,r}^n(L),V_\Z).
$$
It follows from (\ref{calc}) that this group is finitely generated
provided $V\otimes \Q$ does not contain the trivial representation.
Since normal functions are determined by their monodromy, the result
follows.
\end{pf}

If $\V$ contains the trivial representation, the group of normal
functions is an uncountably generated divisible group. For example,
if $\V$ has trivial monodromy, then all such extensions are pulled
back from a point. The set of normal functions is then
$$
\Ext^1_\H(\Z,V_o) \approx JV_o,
$$
where $V_o$ denotes the fiber over the base point.

\begin{theorem}\label{comp}
If, in addition, the fiber over the base point is an irreducible $Sp_g$
module with highest weight $\lambda$ and Hodge weight $m$, then the group
of normal functions $s : \M_{g,r}^n(L) \to J\calV$ is finitely generated
of rank
$$
\dim H^1(\Gamma_{g,r}^n(L), V(\lambda))
=\begin{cases}
1 & \text{if }\lambda = \lambda_3\text{ and } m=-1; \cr
r+n & \text{if }\lambda = \lambda_1\text{ and } m=-1; \cr
0 & \text{ otherwise}.
\end{cases}
$$
\end{theorem}

The upper bounds for the rank of the group of normal functions follow
from (\ref{fin_gen}), (\ref{calc}), and the fact that the  monodromy
representation associated
to a normal function has to be a morphism of variations of mixed Hodge
structure (\ref{morph}). It remains to show that these upper bounds are
achieved. We do this by explicitly constructing normal functions.

Multiples of the generators mod torsion of the normal functions associated
to $V(\lambda_1)$ can be pulled back from $\M_g^1(L)$ along the $n+r$
forgetful maps $\M_{g,r}^n(L) \to \M_g^1(L)$. There the normal function
can be taken to be the one that takes $(C,x)$ to the point
$(2g-2)x -\kappa_C$ of $\Pic^0 C$, where $\kappa_C$ denotes the
canonical class of $C$.

A multiple of the normal function associated to $\lambda_3$ can be pulled
back from $\M_g(L)$ along the forgetful map $\M_{g,r}^n(L) \to \M_g(L)$.
We will describe how this normal function over $\M_g(L)$ arises
geometrically.  If $C$ is a smooth projective curve of genus $g$ and
$x\in C$, we have the Abel-Jacobi mapping
$$
\nu_x : C \to \Jac C.
$$
Denote the algebraic 1-cycle ${\nu_x}_\ast C$ in $\Jac C$ by $C_x$.
Denote the cycle $i_\ast C_x$ by $C_x^-$, where $i : \Jac C \to \Jac C$
takes $u$ to $-u$. The cycle $C_x - C^-$ is homologous to zero, and
therefore defines a point $\etilde(C,x)$ in $JH_3(\Jac C,\Z(-1))$.
Pontrjagin product with the class of $C$ induces a homomorphism
$$
A:\Jac C \to JH_3(\Jac C,\Z(-1)).
$$
Denote the cokernel of $A$ by $JQ(\Jac C)$. It is not difficult to show
that
$$
\etilde(C,x) - \etilde(C,y) = A(x-y).
$$
It follows that the image of $\etilde(C,x)$ in $JQ(\Jac C)$ is
independent of $x$. The image will be denoted by $e(C)$.

The primitive decomposition
$$
H_3(\Jac C,\Q) = H_1(\Jac C,\Q) \oplus PH_3(\Jac C,\Q)
$$
is the decomposition of $H_3(\Jac C)$ into irreducible $Sp_g$
modules; the highest weights of the pieces being $\lambda_1$
and $\lambda_3$, respectively.

Fix a level $L$ so that $\M_g(L)$ is smooth. The union of the
$JQ(\Jac C)$ form the bundle $\J_{\lambda_3}$ of intermediate jacobians
over $\M_g(L)$ associated to the variation of Hodge structure of weight
$-1$ and monodromy the third fundamental representation $V(\lambda_3)$
of $Sp_g$.

\begin{theorem}
The section $e$ of $\J_{\lambda_3}$ is a normal function. Every other
normal function associated to this bundle is, up to torsion, a half
integer multiple of $e$.
\end{theorem}

\begin{pf}
This result is essentially proved in \cite{hain:completions}. We give
a brief sketch.

To see that $e$ is a normal function, consider the
bundle of intermediate jacobians $JH_3(\Jac C,\Z(-1))$ over
$\M_g^1(L)$. It follows from (\ref{norm}) that $(C,x) \mapsto
\etilde(C,x)$ is a normal function. The argument on page 97 of
\cite{hain:completions}, shows that there is a canonical quotient
of the variation corresponding to $\etilde$. (It is the extension
$E$ in display 10 of \cite{hain:completions}.) This variation does
not depend on the base
point $x$, and is therefore constant along the fibers of
$\M_g^1(L) \to \M_g(L)$. It follows that this quotient variation is the
pullback of a variation on $\M_g(L)$. This quotient variation is
classified by $e$. It follows that $e$ is a normal function.

Each normal function $f$ associated to this bundle of intermediate
jacobians induces an $L$ equivariant homomorphism
$$
f_\ast : H_1(T_g,\Z) \to H_1(JQ,\Z) \approx \Lambda^3 H_1(C,\Z)/H_1(C,\Z).
$$
It follows from monodromy computation in
\cite[(4.3.5)]{hain:heights} (see also \cite[(6.3)]{hain:completions})
that $e_\ast$ is twice the Johnson homomorphism.
$$
\tau_g : H_1(T_g,\Z) \to \Lambda^3 H_1(C,\Z)/H_1(C,\Z).
$$
Since this homomorphism is primitive --- i.e., not a non-trivial
integral multiple of another such normal function, all other normal
functions associated to $\lambda_3$ must have monodromy representations
which are half integer multiples of that of $e$. As we have seen in
the proof of (\ref{fin_gen}), such normal functions are determined, up
to torsion, by their monodromy representation. The result follows.
\end{pf}

I don't know how to realize $e/2$ as a normal function in this sense.
But I do know to construct a more general kind of normal function
associated to the 1-cycle $C$ in $\Jac C$ that does realize
$e/2$. It is a section of a bundle whose fiber over $C$ is a principal
$JQ(\Jac C)$ bundle. The details may be found in
\cite[p.~92]{hain:completions}.

\begin{remark}
Using the results in Section \ref{tech} and Theorem \ref{comp}, one
can easily show that the rank of the group of normal
functions in the theorem above is
$$
\dim \Gamma\Hom_{Sp_g(\Q)}(H_1(T_{g,r}^n,\Q),V_{\Q,C}),
$$
where $H_1(T_{g,r}^n)$ is given the Hodge structure of weight $-1$
described in \S \ref{tech}.
\end{remark}

\section{Technical Results on Variations over $\M_g$}
\label{tech}

In this section, we prove several technical facts about variations of
mixed Hodge structure over moduli spaces of curves that were used
in Section \ref{normal_functions}. Throughout we will assume that $L$ has
been chosen so that $\Gamma_{g,r}^n(L)$ is torsion free.

\begin{proposition}\label{exist-unique}
The local system $\V(\lambda)$ over $\M_{g,r}^n(L)$ associated
to the irreducible representation of $Sp_g$ with highest weight
$\lambda$ underlies a good $\Q$ variation of (mixed) Hodge structure,
and this variation is unique up to Tate twist.
\end{proposition}

\begin{pf}
First observe that the local system $\bH$ corresponding to the
fundamental representation $V(\lambda_1)$ occurs as a variation of Hodge
structure over $\M_{g,r}^n(L)$ of weight 1; it is simply the local system
$R^1\pi_\ast\Q$ associated to the universal curve $\calC \to \M_{g,r}^n(L)$.
The existence of the structure of a good variation of Hodge structure
on the local system corresponding the the $Sp_g$ module with
highest weight $\lambda$ now follows using Weyl's construction of the
irreducible representations of $Sp_g$---see, for example,
\cite[\S 17.3]{fulton-harris}.

To prove uniqueness, suppose that $\V$ and $\V'$ are both good variations
of mixed Hodge structure corresponding to the same irreducible $Sp_g$ module.
{}From Saito \cite{saito}, we know that
$$
\Hom_{\Gamma_{g,r}^n(L)}(\V,\V')
$$
has a mixed Hodge structure. By Schur's lemma, this group is
one dimensional. It follows that this group is isomorphic to $\Q(n)$
for some $n$. It follows that $\V' = \V(n)$.
\end{pf}

\begin{proposition}\label{decomp}
If $\E$ is a good variation of $\Q$ mixed Hodge structure over
$\M_g(L)$ whose monodromy representation factors through a rational
representation of the algebraic group $Sp_g$, then for each dominant
integral weight $\lambda$ of $Sp_g$, the $\lambda$-isotypical part
$\E_\lambda$ of $\E$ is a good variation of mixed Hodge structure.
Consequently,
$$
\E = \bigoplus_\lambda \E_\lambda
$$
in the category of good variations of $\Q$ mixed Hodge structure
over $\M_g(L)$. Moreover, for each $\lambda$, there is a mixed Hodge
structure $A_\lambda$ such that
$\E_\lambda = A_\lambda \otimes \V(\lambda)$.
\end{proposition}

\begin{pf}
Fix $\lambda$, and let $\V(\lambda) \to \M_g(L)$ be a variation of
Hodge structure whose fiber over some fixed base point is the
irreducible $Sp_g$ module with highest weight $\lambda$.
It follows from Saito's work \cite{saito} that
$$
A_\lambda := \Hom_{\Gamma_g(L)}\big(\V(\lambda), \E\big) =
H^0\big(\M_g(L),\Hom_\Q(\V(\lambda),\E)\big)
$$
is a mixed Hodge structure.  Let
$$
\E' = \bigoplus_\lambda A_\lambda \otimes \V(\lambda).
$$
This is a good variation of mixed Hodge structure which is isomorphic
to $\E$ as a $\Q$ local system. Now
$$
\Hom_{\Gamma_g(L)}(\E',\E) = \bigoplus_\lambda A_\lambda^\ast
\otimes \Hom_{\Gamma_g(L)}(\V(\lambda),\E) = \bigoplus_\lambda
\Hom_\Q(A_\lambda,A_\lambda).
$$
The element of this group which corresponds to
$\id : A_\lambda \to A_\lambda$ in each factor is an isomorphism
of local systems and an element of
$$
\Gamma\Hom_{\Gamma_g(L)}(\E',\E).
$$
It is therefore an isomorphism of variations of mixed Hodge structure.
\end{pf}

The local system
$$
\left\{H_1(T_{g,r}^n)\right\}
$$
over $\M_{g,r}^n(L)$ naturally underlies a variation of mixed Hodge
structure of weight $-1$. The $\lambda_1$ isotypical component
is simply $r+n$ copies of the variation $\V(\lambda_1)$. We shall
denote this variation by $\bH_1(T_{g,r}^n)$.

\begin{proposition}\label{morph}
Suppose that $\V$ is a variation of mixed Hodge structure over
$\M_{g,r}^n(L)$ whose monodromy representation factors through a
rational representation of $Sp_g$. If $\E$ is an extension of
$\Q$ by $\V$ in the category of variations of mixed Hodge structure
over $\M_{g,r}^n(L)$, then the restriction of the monodromy
representation to $H_1(T_{g,r}^n)$,
$$
\bH_1(T_{g,r}^n) \to \V,
$$
is a morphism of variations of mixed Hodge structure.
\end{proposition}

\begin{pf}
It suffices to prove the assertion for $\Q$ variations of mixed
Hodge structure. We will prove the case where $n=r=0$; the proofs
of the other cases being similar.

If the monodromy representation of $\E$ is trivial, the
result is trivially true. So we shall assume that the monodromy
representation is non-trivial.

Using the previous result, we can write
$$
\V = \bigoplus_\lambda \V_\lambda
$$
as variations of mixed Hodge structure over $\V$. By pushing
out the extension
$$
0 \to \V \to \E \to \Q \to 0
$$
along the projection $\V \to \V_{\lambda_3}$ onto the $\lambda_3$
isotypical component, we obtain an extension
$$
0 \to \V_{\lambda_3} \to \E' \to \Q \to 0.
$$
It follows from Johnson's computation that the restricted monodromy
representation of $\E$ factors through that of $\E'$:
$$
\bH_1(T_g) \to \V_{\lambda_3} \to \V.
$$
We have therefore reduced to the case where $\V=\V_{\lambda_3}$.

Let $\V(\lambda_3)$ be the unique variation of Hodge structure of
weight $-1$ over $\M_g(L)$ with monodromy representation given by
$\lambda_3$. Let $\bS$ be the variation of mixed Hodge structure over
$\M_g(L)$ given by the cycle $C - C^-$ that was constructed in Section
\ref{normal_functions}. It is an extension of $\Q$ by $\V(\lambda_3)$.

By \cite{saito}, the exact sequence
$$
0 \to \Hom_{\Gamma_g(L)}(\bS,\V_{\lambda_3}) \to
\Hom_{\Gamma_g(L)}(\bS,\E') \to \Hom_{\Gamma_g(L)}(\bS,\Q)
$$
is a sequence of mixed Hodge structures. The most right hand group is
easily seen to be isomorphic to $\Q(0)$; it is generated by the
projection $\bS \to \Q$. The left hand group is easily seen to be
zero. It follows that
$$
\Hom_{\Gamma_g(L)}(\bS,\E') \approx \Q(0).
$$
Since the monodromy representation of $\bS$ is a morphism, it follows
that the monodromy representations of $\E'$ and $\E$ are too.
\end{pf}

\section{The Harris-Pulte Theorem}

As an application of the classification of normal functions above, we
give a new proof of the Harris-Pulte theorem which relates the mixed
Hodge structure on $\pi_1(C,x)$ to the normal function of the cycle
$C_x-C_x^-$ when $g\ge 3$.  The result we obtain is slightly stronger.

Fix a level so that $\Gamma_g^1(L)$ is torsion free.
Denote by $\bL$ the $\Z$ variation of Hodge structure of weight
$-1$ over $\M_g^1(L)$ whose fiber over the pointed curve $(C,x)$
is $H_1(C)$. Denote the corresponding holomorphic vector bundle
by $\L$. The cycle $C_x - C_x^-$ defines a normal function $\zeta$
which is a section of
$$
J\Lambda^3\L \to \M_g^1(L).
$$

Denote the integral group ring of $\pi_1(C,x)$ by $\Z\pi_1(C,x)$,
and its augmentation ideal by $I(C,x)$, or $I$ when there is no
possibility of confusion. There is a canonical mixed
Hodge structure on the truncated augmentation ideal
$$
I(C,x)/I^3.
$$
(See, for example, \cite{hain:geom}.)
It is an extension
$$
0 \to H_1(C)^{\otimes 2}/q \to I(C,x)/I^3 \to H_1(C) \to 0,
$$
where $q$ denotes the symplectic form. Tensoring with $H_1(C)$
and pulling back the resulting extension along the map
$\Z \to H_1(C)^{\otimes 2}$, we obtain an extension
$$
0 \to H_1(C)\otimes \left(H_1(C)^{\otimes 2}/q\right)
\to E(C,x) \to \Z \to 0.
$$
Since the set of $I(C,x)$ form a good variation of mixed Hodge
structure over $\M_g^1(L)$ (\cite{hain:dht}), the set of $E(C,x)$
form a good variation of mixed Hodge structure $\E$ over $\M_g^1(L)$.
It therefore determines a normal function $\rho$ which is a section of
$$
J\L\otimes\left(\L^{\otimes 2}/q\right) \to \M_g^1(L).
$$

Define the map
$$
\Phi: J\Lambda^3 \L \to J\L\otimes(\L^{\otimes 2}/q)
$$
to be the one induced by the map
$$
\Lambda^3 \bL \to \bL^{\otimes 3} \to
\bL \to \bL\otimes(\bL^{\otimes 2}/q);
$$
the first map is defined by
$$
x_1 \wedge x_2 \wedge x_3 \mapsto \sum_\sigma
\sgn(\sigma) x_{\sigma(1)} \otimes x_{\sigma(2)} \otimes x_{\sigma(3)}
$$
where $\sigma$ ranges over all permutations of $\{1,2,3\}$.

Our version of the Harris-Pulte Theorem is:

\begin{theorem}
The image of $\zeta$ under $\Phi$ is $2\rho$.
\end{theorem}

\begin{pf}
The proof uses (\ref{rigidity}). It is a straightforward consequence
of (\ref{equality}) that the monodromy representations of $\Phi(\zeta)$
and $2\rho$ are equal. It is also a straightforward matter
to use functoriality to show that both $\Phi(\zeta)$ and $2\rho$ vanish
at $(C,x)$ when $C$ is hyperelliptic and $x$ is a Weierstrass point
(cf.\ \cite[(7.5)]{hain:geom}.)
\end{pf}

\section{The Franchetta Conjecture for Curves with a Level}
\label{franchetta}

Suppose that $L$ is a finite index subgroup of $Sp_g(\Z)$, not necessarily
torsion free. Denote the generic point of $\M_g(L)$ by $\eta$. There is a
universal curve defined generically over $\M_g(L)$. Denote its fiber over
$\eta$ by $\calC_g(L)_\eta$. In the statement below, $S$ denotes a
compact oriented surface of genus $g$.

\begin{theorem}\label{franch_conj}
For all $g\ge 3$ and all finite index subgroups $L$ of $Sp_g(\Z)$, the
group $\Pic \calC_g(L)_\eta$ is finitely generated of rank 1. The
torsion subgroup is isomorphic to $H^0(L,H_1(S,\Q/\Z))$. Modulo torsion,
either it is generated by the canonical bundle, or by a divisor of degree
$g-1$.
\end{theorem}

This has a concrete statement when $L=Sp_g(\Z)(l)$, the congruence
subgroup of level $l$ of $Sp_g(\Z)$.  It is not difficult to show
that the only torsion points of $\Jac S$ invariant under $L$ are
the points of order $l$. That is,
$$
H^0(L,H_1(S,\Q/\Z)) \approx H_1(S,\Z/l\Z).
$$
In this case we shall denote $\calC_g(L)_\eta$ by $\calC_g(l)_\eta$.
During the proof of the theorem, we will show that, mod torsion,
$\Pic \calC_g(l)_\eta$ is generated by a theta characteristic when $l$
is even. Combining this with the theorem, we have:

\begin{corollary}
If $g\ge 3$, then for all $l\ge 0$, $\Pic \calC_g(l)_\eta$ is a
finitely generated group of rank one with torsion subgroup
isomorphic to $H_1(S,\Z/l\Z)$. Modulo torsion, $\Pic \calC_g(l)_\eta$
is generated by a theta  characteristic when $l$ is even, and by the
canonical bundle when $l$ is odd. \qed
\end{corollary}

The case $g=2$, if true, should follow from Mess's computation of
$H_1(T_2)$ \cite{mess}. One should note that Mess proved that $T_2$
is a countably generated free group.

\begin{pf*}{Sketch of proof of Theorem \ref{franch_conj}}
We first suppose that $L$ is torsion free. In this case, the universal
curve is defined over all of $\M_g(L)$.  Denote the restriction of
it to a Zariski open subset $U$ of $\M_g(L)$ by $\calC_g(L)_U$. Set
$$
\Pic_{\rel U} \calC_g(L) = \coker\{\Pic U \to \Pic \calC_g(U)\}.
$$
Then
$$
\Pic \calC_g(L)_\eta =
\lim_{\stackrel{\longrightarrow}{U}}
\Pic_{\rel U} \calC_g(L),
$$
where $U$ ranges over all Zariski open subsets of $\M_g(L)$.
There is a natural homomorphism
$$
\mathrm{deg}: \Pic_{\rel U}\calC_g \to \Z
$$
given by taking the degree on a fiber. Denote $\mathrm{deg}^{-1}(d)$
by $\Pic_{\rel U}^d \calC_g(L)$.

We first compute $\Pic^0 \calC_g(L)_\eta$. Each element of this group
can be represented by a line bundle over $\calC_g(L)_U$ whose
restriction to each fiber of $\pi : \calC_g(L)_U \to U$ is topologically
trivial. This line bundle has a section. By tensoring it with the
pullback of a line bundle on $U$, if necessary, we may assume that the
divisor of this section intersects each fiber of $\pi$ in only
a finite number of points. We therefore obtain a normal function
$$
s : U \to \Pic_{\rel U}^0 \calC_g(L).
$$
Since the associated variation of Hodge structure is the unique one
of weight $-1$ associated to $V(\lambda_1)$, it follows from
(\ref{comp}) and (\ref{extend}) that this normal function is torsion.
It follows that
$$
\Pic^0 \calC_g(L)_\eta = \Pic_{\rel U}^0 \calC_g(L)
= H^0(L,H_1(S,\Q/\Z)).
$$
Since this group is isomorphic to $H_1(S,\Z/l\Z)$ when $L$ is the
congruence $l$ subgroup of $Sp_g(\Z)$, and since every finite index
subgroup of $Sp_g(\Z)$ contains a congruence subgroup by \cite{b-m-s},
it follows that $\Pic^0 \calC_g(L)_\eta$ is finite for all $L$.

The relative dualizing sheaf $\omega$ of $\calC_g(L)_U$ gives an
element of $\Pic^{2g-2} \calC_g(L)_\eta$. Denote the greatest common
divisor of the degrees of elements of $\Pic \calC_g(L)_\eta$ by $d$.
Observe that $d$ divides $2g-2$. Let $m=(2g-2)/d$. We will show that
$m=1$ or 2.

Choose an element $\delta$ of $\Pic^d \calC_g(L)_\eta$. Then
$$
\omega - m\delta \in \Pic^0 \calC_g(L)_\eta
$$
and is therefore torsion of order $k$, say. Replace $L$ by
$$
L' = L \cap Sp_g(\Z)(km).
$$
Observe that the natural map
$$
\Pic^0 \calC_g(L)_\eta \to \Pic^0 \calC_g(L')_\eta
$$
is injective. We can find
$$
\mu \in \Pic^0 \calC_g(L')_\eta
$$
such that $m\mu = \omega - m\delta$. Then $\delta + \mu$ is an $m$th
root of the canonical bundle $\omega$. It appears to be well known
that the only non-trivial roots of the canonical bundle that can be
defined over $\M_g(L)$ are square roots. In any case, this follows from
the result (\ref{action}) in the next section. This implies that
$m$ divides 2, as claimed.

It follows from (\ref{action})) that square roots of the canonical
bundle are defined over $\M_g(l)$ if and only if $l$ is even. Combined
with the argument above, this shows that, mod torsion,
$\Pic^0 \calC_g(l)_\eta$ is  generated by $\omega$ if $l$ is odd, and by
a square root of $\omega$ if $l$ is even.

Our final task is to reduce the general case to that where $L$ is
torsion free. For arbitrary $L$, we have
$$
\Pic \calC_g(L)_\eta =
\lim_{\stackrel{\longrightarrow}{U}}
\Pic_{\rel U} \calC_g(L),
$$
where $U$ ranges over all smooth Zariski open subsets of $\M_g(L)$.
Choose a torsion free finite index normal subgroup $L'$ of $L$ and a
smooth Zariski open subset $U$ of $\M_g(L)$. Denote the inverse image
of $U$ in $\M_g(L')$ by $U'$. Then the projection $U' \to U$ is a
Galois cover with Galois group $G=L/L'$. It follows that
$$
\Pic_{\rel U} \calC_g(L) = \Pic_{\rel U'} \calC_g(L')^G.
$$
Since $\pi_1(U)$ surjects onto $\Gamma_g(L)$, and therefore onto $L$,
the result follows.
\end{pf*}

Denote the universal curve over the generic point $\eta$ of
$\M_{g,r}^n(l)$ by $\calC_{g,r}^n(l)_\eta$. The proof of the following
more general result is  similar to that of Theorem \ref{franch_conj}.
\begin{theorem}
If $g\ge 3$, then for all $l\ge 0$, $\Pic\calC_{g,r}^n(l)_\eta$ is a
finitely generated group of rank $r+n+1$ whose torsion subgroup
isomorphic to $H_1(S,\Z/l\Z)$. Each of the $n$ marked points and the
anchor point of each of the $r$ marked cotangent vectors gives an element of
$\Pic^1\calC_{g,r}^n(l)_\eta$. The pairwise differences of these points
generate a subgroup of $\Pic^0\calC_{g,r}^n(l)_\eta$ of rank $r+n-1$.
Moreover, $\Pic^0\calC_{g,r}^n(l)_\eta$ is generated by these differences
modulo torsion. Modulo $\Pic^0\calC_{g,r}^n(l)_\eta$,
$\Pic\calC_{g,r}^n(l)_\eta$ is generated by the class of one of the
distinguished points together with a theta characteristic when $l$ is
even, and by the canonical divisor when $l$ is odd. \qed
\end{theorem}

Note that the independence of the pairwise difference of the points
follows from the discussion following Theorem \ref{comp}.

\section{The Monodromy of Roots of the Canonical Bundle}
\label{root}

In this section we compute the action of $\Gamma_g$ on the set of
$n$th roots of the canonical bundle of a curve of genus $g$.
This action has also been computed by P.~Sipe \cite{sipe}, but
in quite a different form.

If $L$ is an $n$th root of the tangent bundle of a smooth projective
curve $C$, then its dual is an $n$th root of the canonical bundle.
That is, there is a one-one correspondence between $n$th roots of
the canonical bundle and $n$th roots of the tangent bundle of a curve.
As it is more convenient, we shall work with roots of the tangent bundle.

The first point is that roots of the tangent bundle are determined
topologically (cf.\ \cite[\S 3]{atiyah} and \cite{sipe}): denote the
$\C^\ast$ bundle associated to the holomorphic tangent bundle $TC$ of
$C$ by $T^\ast$. Indeed, an $n$th root of $TC$ is a cyclic covering of
$T^\ast$ of degree $n$ which has degree $n$ on each fiber. The complex
structure on such a covering is uniquely determined by that on $T^\ast$.

The first Chern class of $TC$ is $2-2g$. So if $R$ is an $n$th root of
$K$, we have that $n$ divides $2g-2$. Since the Euler class of $T^\ast$
is $2-2g$, it follows from the Gysin sequence that there is a short exact
sequence
\begin{equation}\label{extension}
0 \to \Z/n\Z \to H_1(T^\ast,\Z/n\Z) \to H_1(C,\Z/n\Z) \to 0
\end{equation}
By covering space theory, an $n$th root of $TC$ is determined by
a homomorphism
$$
H_1(T^\ast,\Z/n\Z) \to \Z/n\Z
$$
whose composition with the inclusion
$\Z/n\Z \hookrightarrow H_1(T^\ast,\Z/n\Z)$ is the identity. That is,
we have the following result:

\begin{proposition}
There is a natural one-to-one correspondence between $n$th roots of the
the canonical bundle of $C$ and splittings of the sequence
(\ref{extension}). \qed
\end{proposition}

Throughout this section, we will assume $g \ge 3$.
Denote the set of $n$th roots of $TC$ by $\Theta_n$.  This is a
principal $H_1(C,\Z/n\Z)$ space. The automorphisms of this affine
space is an extension
$$
0 \to H_1(C,\Z/n\Z) \to \Aut \Theta_n
\stackrel{\pi}{\to} GL_{2g}(\Z/n\Z) \to 1;
$$
the kernel being the group of translations by elements of $H_1(C,\Z/n\Z)$.
The mapping class group acts on $\Theta_n$, so we have a homomorphism
$$
\Gamma_g \to \Aut \Theta_n.
$$
The composite of this homomorphism with $\pi$ is the
reduction mod $n$
$$
\rho_n : \Gamma_g \to Sp_g(\Z/n\Z)
$$
of the natural homomorphism. Denote the subgroup
$\pi^{-1}(Sp_g(\Z/n\Z))$ of $\Aut \Theta_n$ by $\calK_n$. It follows
that the action of $\Gamma_g$ on $\Theta_n$ factors through a
homomorphism
$$
\theta_n : \Gamma_g \to \calK_n
$$
whose composition with the natural projection
$\calK_n \to Sp_g(\Z/n\Z)$ is $\rho_n$. In order to determine
$\theta_n$, we will need to compute its restriction
$$
\theta_n : H_1(T_g) \to H_1(C,\Z/n\Z)
$$
to the Torelli group. First some algebra.

\begin{proposition}
There is a natural homomorphism
$$
\psi_g : H_1(T_g,\Z) \to H_1(C,\Z/(g-1)\Z).
$$
\end{proposition}

\begin{pf}
By (\ref{tau_g}), there is a natural homomorphism
$$
\tau_g : H_1(T_g,\Z) \to
\Lambda^3 H_1(C,\Z)/\big([C]\times H_1(C,\Z)\big).
$$
Here we view $\Lambda^\dot H_1(C)$ as the homology of $\Jac C$ and
$[C]$ denotes the homology class of the image of $C$ under
the Abel-Jacobi map.  There is also a natural homomorphism
$$
p : \Lambda^3 H_1(C,\Z) \to H_1(C,\Z)
$$
defined by
$$
p : x\wedge y \wedge z \mapsto
(x\cdot y)\, z + (y\cdot z)\, x + (z\cdot x)\, y.
$$
It is easy to see that the composite
$$
H_1(C,\Z) \stackrel{[C]\times}{\longrightarrow}
\Lambda^3 H_1(C,\Z) \stackrel{p}{\to} H_1(C,\Z)
$$
is multiplication by $g-1$. It follows that $p$ induces a homomorphism
$$
\pbar : \Lambda^3 H_1(C,\Z) {\to} H_1(C,\Z/(g-1)\Z).
$$
The homomorphism $\psi_g$ is the composite $\pbar\circ\tau_g$.
\end{pf}

Call a translation of $\Theta_n$ {\it even\/} if it is translation
by an element of $2H^1(C,\Z/n\Z)$. If $n$ is odd, this is the set
of all translations. If $n=2m$, this is the proper subgroup of
$H^1(C,\Z/n\Z)$ isomorphic to $H^1(C,\Z/m\Z)$. It is not difficult to
see that there is a unique subgroup of $\calK_n$ that is an extension
of $Sp_g(\Z/n\Z)$ by the even translations. We shall denote it by
$\calK_n^{(2)}$.

\begin{theorem}\label{action}
The image of the natural homomorphism $ \theta_n : \Gamma_g \to \calK_n$
is $\calK_n^{(2)}$. The restriction of $\theta_n$ to $T_g$ is the
composite of $\psi_g$ with the homomorphism
$$
H_1(C,\Z/(g-1)\Z) @>r>> H_1(C,\Z/n\Z) @>{PD}>> H^1(C,\Z/n\Z),
$$
where $r(k)$ equals $2k$ mod $n$ and $PD$ denotes Poincar\'e duality.
In particular, the Torelli group acts trivially on $\Theta_n$ if and
only if $n$ divides 2.
\end{theorem}

\begin{pf}
First, it was proven by Johnson in \cite{johnson_2} that the kernel
of the composite
$$
 T_g \to H_1(T_g) \stackrel{\tau_g}{\to}
\Lambda^3 H_1(C,\Z)/\big([C]\times H_1(C,\Z)\big)
$$
is generated by Dehn twists on separating simple closed curves. Using
this, it is easy to check that the restriction of $\theta_n$ to $T_g$
factors through $\tau_g$. In \cite{johnson_3}, Johnson shows that
$T_g$ is generated by Dehn twists on a bounding pair of disjoint
simple closed curves.\footnote{Actually, all we need is that
$\Lambda^3 H_1(C,\Z)/\big([C]\times H_1(C,\Z)\big)$ be generated
by the images under $\tau_g$ by such bounding pair maps.
This is easily checked directly.}

Now suppose that $\varphi$ is such a bounding pair map. There are
two disjoint imbedded circles $A$ and $B$ such that $\varphi$ equals
a positive Dehn twist about $A$ and a negative one about $B$.
When we cut $C$ along $A\cup B$, we obtain two surfaces, of
genera $g'$ and $g''$, say. Choose one of these components, and
let $a$ be the cycle obtained by orienting $A$ so that it is a
boundary component of this component. It is not difficult to show
that the image of $\varphi$ under $\psi_g$ equals
$$
-g'\, [a] \in H_1(C,\Z/(g-1)\Z).
$$
where $g'$ is the genus of the chosen component. Since $g' + g'' = g-1$,
this is well defined. Next, one can use Morse theory to show that the
image of this same bounding pair map in $H^1(C,\Z/n\Z)$ is $-2 g'PD(a)$,
from which the result follows.  Full details of these computation will
appear elsewhere.
\end{pf}

\begin{corollary}
The only roots of the canonical bundle defined over Torelli space
are the canonical bundle itself and its $2^{2g}$ square roots. \qed
\end{corollary}

\begin{remark}
The homomorphism $\theta_{2g-2} :\Gamma_g \to \calK_{2g-2}$ appears
in Morita's work --- cf. \cite[\S 4.A]{morita}.
\end{remark}

\section{Heights of Cycles defined over $\M_g(L)$}
\label{heights}

Suppose that $X$ is a compact K\"ahler manifold of dimension $n$
and that $Z$ and $W$ are two homologically trivial algebraic cycles
in $X$ of dimensions $d$ and $e$, respectively. Suppose that
$d+e = n-1$ and that $Z$ and $W$ have disjoint supports. Denote the
current associated to $W$ by $\delta_W$. It follows from the
$\del\delbar$-Lemma that there is a current $\eta_W$ of type
$(d,d)$ that is smooth away from the support of $Z$ and
satisfies
$$
\del\delbar\eta_W = \pi i \delta_W.
$$
The (archimedean) height pairing between $Z$ and $W$ is defined by
$$
\langle Z,W \rangle = - \int_Z \eta_W.
$$
This is a real-valued, symmetric bilinear pairing on such disjoint
homologically trivial cycles. It is important in number theory
(cf.\ \cite{beilinson:height}).

Now suppose that
$$
X \to \M_g(L)
$$
is a family of smooth projective varieties of relative dimension $n$.
Suppose that $Z\to \M_g(L)$ and $W\to \M_g(L)$ are families of
algebraic cycles in $X$ of relative dimensions $d$ and  $e$, respectively,
where $d+e = n-1$. Denote the fiber of $X$, $Z$ and $W$ over $C \in \M_g(L)$
by $X_C$, $Z_C$ and $W_C$, respectively.  Suppose that $Z_C$ and $W_C$ are
homologically trivial in $X_C$ and that they have disjoint supports for
generic $C\in \M_g(L)$.

We shall suppose that $L$ has been chosen so that every curve has
two distinguished theta characteristics $\alpha$ and $\alpha+\delta$,
where $\delta$ is a non-zero point of order 2 in $\Jac C$. We shall
also suppose that $g$ is odd and $\ge 3$. Write $g$ in the form $g=2d+1$.

Denote the difference divisor
$$
\left\{x_1+\dots + x_d -y_1 -\dots - y_d : x_j,y_j\in  C\right\}
$$
in $\Jac C$ by $\Delta$, and the theta divisor
$$
\left\{x_1+\dots +x_{2d} - \alpha : x_j \in C\right\}
$$
in $\Jac C$ by $\Theta_\alpha$. By \cite[(4,1.2)]{hain:heights}, there
is a rational function $f_C$ on $\Jac C$ whose divisor is
$$
\Delta - {2d \choose d} \Theta_\alpha.
$$
Denote the unique invariant measure of total mass one on $\Jac C$
by $\mu$.

\begin{theorem}\label{height}
Suppose that $g$ is odd and $\ge 3$. Suppose that $Z$ and $W$ are
families of homologically trivial cycles over $\M_g(L)$ in a family of
smooth projective varieties $p:X\to \M_g(L)$, as above. If the
monodromy of the local system $R^{2d+1}p_\ast \Q_X$ factors through
a rational representation of $Sp_g$, then there is a rational function
$h$ on $\M_g(L)$, and rational numbers $a$ and $b$ such that
$$
\langle Z_C,W_C\rangle = a \left(
\log |h(C)| + 2b\left(\log |f_C(\delta)|
- \int_{\Jac C}\log |f_C(x)| d\mu(x)\right)\right).
$$
\end{theorem}

The numbers $a$ and $b$ are topologically determined, as will become
apparent in the proof.  The divisor of $h$ is computable when one has
a good understanding of how the cycles $Z$ and $W$ intersect.
One should be able to derive a similar formula for even $g$ using Bost's
general computation of the height in \cite{bost} and results from
\cite{hain:completions}.

The proof of Theorem \ref{height} occupies the remainder of this
section. We only give a sketch. We commence by defining two
algebraic cycles in $\Pic^d C$. For $D\in \Jac C$, let $C^{(d)}_D$ be
the $d$-cycle in $\Pic^d C$ obtained by pushing forward the
fundamental class of the $d$th symmetric power of $C$ along the
map
$$
\{x_1,\dots,x_d\} \mapsto x_1 + \dots + x_d + D.
$$
Let $i$ be the automorphism of $\Pic^d C$ defined by
$i : x \mapsto \alpha - x$.
Define
$$
Z_D = C^{(d)}_D - i_\ast C^{(d)}_D.
$$
This is a homologically trivial $d$-cycle in $\Pic^d C$.

{}From \cite{bost} and \cite{hain:heights}, we know that
$$
\langle Z_0 , Z_\delta \rangle =
2\log |f_C(\delta)| - 2\int_{\Jac C}\log |f_C(x)| d\mu(x).
$$
So the content of the theorem is that there is a rational function
$h$ on $\M_g(L)$ and rational numbers $a$ and $b$ such that
$$
\langle Z,W \rangle =
a\left(\log |h(C)| + b\, \langle Z_0 , Z_\delta \rangle\right).
$$
The basic point, as we shall see, is that, up to torsion, all normal
functions over $\M_g(L)$ are half integer multiples of that of $C-C^-$,
as was proved in Section \ref{normal_functions}.

We will henceforth assume that the reader is familiar with the content
of \cite[\S 3]{hain:heights}. We will briefly review the most relevant
points of that section.

A {\it biextension} is a mixed Hodge structure $B$ with only three
non-trivial weight graded quotients: $\Z$, $H$, and $\Z(1)$, where
$H$ is a Hodge structure of weight $-1$. The isomorphisms
$$
Gr^W_{-2}B\approx \Z(1)\text{ and }Gr^W_0B\approx \Z
$$
are considered to be part of the data of the biextension. If one replaces
$\Z$ by $\R$ in this definition, one obtains the definition of a real
biextension. To a biextension $B$ one can associate a real number
$\nu(B)$, called the {\it height} of $B$. It depends only on the
associated real biextension $B\otimes\R$.

To a pair of disjoint homologically trivial cycles in a
smooth projective variety $X$ satisfying
$$
\dim Z + \dim W + 1 = \dim X,
$$
there is a canonical biextension $B_\Z(Z,W)$, whose
weight graded quotients are
$$
\Z,\quad H_{2d+1}(X,\Z(-d)),\quad \Z(1),
$$
where $d$ is the dimension of $Z$.
The extensions
$$
0 \to H_{2d+1}(X,\Z(-d)) \to B_\Z(Z,W)/\Z(1) \to \Z \to 0
$$
and
$$
0 \to \Z(1) \to W_{-1}B_\Z(Z,W) \to H_{2d+1}(X,\Z(-d)) \to 0
$$
are the those determined by $Z$ (directly), and $W$
(via duality) \cite[(3.3.2)]{hain:heights}.  We have
$$
\nu(B_\Z(Z,W)) = \langle Z,W \rangle.
$$

The first step in the proof is to reduce the size of the biextension.
Suppose that $\Lambda=\Z$ or $\R$, and that $B$ is a
$\Lambda$-biextension with weight $-1$ graded quotient $H$.
Suppose that there is an inclusion $i: A\hookrightarrow H$ of
$\Lambda$ mixed Hodge structures. Pulling back the extension
$$
0 \to \Lambda(1) \to W_1B \to H \to 0
$$
along $i$, we obtain an extension
$$
0 \to \Lambda(1) \to E \to C \to 0.
$$
If this extension splits, there is a canonical lift
$\itilde : C \to B$ of $i$. The quotient $B/C$ is also a $\Lambda$
biextension.

\begin{proposition}\label{prune}
The biextensions $B_\Lambda(Z,W)$ and $B_\Lambda(Z,W)/C$ have the
same height.
\end{proposition}

\begin{pf}
This is a special case of \cite[(5.3.8)]{lear}. It
follows directly from \cite[(3.2.11)]{hain:heights}.
\end{pf}

We will combine this with (\ref{comp}) to prune the biextension
$B_\Z(Z_C,W_C)$ until its weight $-1$ graded quotient is either
trivial or else one copy of $V(\lambda_3)$.

First observe that if $B$ is a biextension and $B'$ a mixed Hodge
substructure of $B$ of finite index, then $B'$ is a biextension and
there is a non-zero integer $m$ such that $\nu(B') = m\nu(B)$. This can be
proved using \cite[(3.2.11)]{hain:heights}.

To prune the biextension $\B(Z,W)$ over $\M_g(L)$, we consider the portion
of the monodromy representation
$$
H_1(T_g) \to \Hom_\Z(Gr^W_{-1}B(Z_C,W_C),\Z(1))
$$
associated to the variation $W_{-1}\B(Z,W)$ over $\M_g(L)$. This map is
$Sp_g$ equivariant. Denote $Gr^W_{-1}\B(Z,W)$ by $\bH$. This monodromy
representation corresponds to the map
$$
\bH \to \left\{H^1(T_g,\Z(1))\right\}
$$
of local systems over $\M_g(L)$ which takes $h\in H_C$ to the
functional $\{\phi \mapsto \phi(h)\}$ on $H_1(T_g)$. For each
$C\in\M_g(L)$ this is a morphism of Hodge structures by (\ref{morph}).
Denote its kernel by $K_C$. These form a variation of Hodge structure
$\K$ over $\M_g(L)$. If the monodromy representation is trivial on
$T_g$, then $\K = \bH$. Otherwise, Schur's lemma implies that $\bH/\K$
is isomorphic $\V(\lambda_3)$ placed in weight $-1$.

We can pull back the extension
$$
0 \to \Q(1) \to W_{-1}\B(Z,W) \to \bH \to 0
$$
along the inclusion $\K\hookrightarrow \bH$ to obtain an extension
\begin{equation}\label{putative}
0 \to \Q(1) \to \E \to \K \to 0.
\end{equation}
If this extension splits over $\Q$, then, by replacing the lattice in
$\B_\Z(Z,W)$ by a commensurable one, we may assume that the
splitting is defined over $\Z$.  This has the effect of multiplying the
height by a non-zero rational number. Once we have done this, the
inclusion $\K \hookrightarrow Gr^W_{-1}\B(Z,W)$ lifts to an inclusion
$\K \hookrightarrow \B(Z,W)$. Using (\ref{prune}), we can replace
$B_\Z(Z_C,W_C)$ by $\B' = \B(Z,W)/\K$ without changing height of
the biextension. For the time being, we shall assume that (\ref{putative})
splits over $\Q$. This is the case, for example, when $\bH$
contains no copies of the trivial representation, as follows from
(\ref{decomp}) since $\K$ is a trivial $T_g$ module by construction.

The weight $-1$ graded quotient of $\B'$
is either trivial or isomorphic to $\V(\lambda_3)$. This
biextension is defined over the open subset $U$ of $\M_g(L)$ where
$Z_C$ and $W_C$ are disjoint. The related variations $W_{-1}\B'$
and $\B'/\Z(1)$ are defined over all of $\M_g(L)$.

If $\K=\bH$, then $\B'$ is an extension of $\Z$ by $\Z(1)$. It therefore
corresponds to a rational function $h$ on $\M_g(L)$ which is defined
on $U$ (cf. \cite[(9.3)]{hain:poly}). It follows from
\cite[(3.2.11)]{hain:heights} that the height of this biextension $\B'$
is $C\mapsto \log|h(C)|$. This completes the proof of the theorem in
this case.

Dually, when the extension
$$
0 \to Gr^W_{-1}\B' \to \B'/\Z(1) \to \Z \to 0
$$
has finite monodromy, there is a rational function $h$ on $\M_g(L)$
such that the height of $B'$, and therefore $B(Z_C,W_C)$, is rational
multiple of $\log|h(C)|$.

We have therefore reduced to the case where $\B'$ has weight graded
$-1$ quotient $\V(\lambda_3)$ and where neither of the extensions
$$
0 \to \V(\lambda_3) \to \B'/\Z(1) \to \Z \to 0
$$
or
$$
0 \to \Z(1) \to W_{-1}\B' \to \V(\lambda_3) \to 0
$$
is torsion. We also have the biextension $\B''$ associated to the cycles
$Z_0$ and $Z_\delta$. It has these same properties. After replacing the
lattices in each by lattices of finite index, we may assume that the
extensions of variations $W_{-1}\B'$ and $W_{-1}\B''$ are isomorphic,
and that the $\B'/\Z(1)$ and $\B''/\Z(1)$ are isomorphic. As in
\cite[(3.4)]{hain:heights}, the
biextensions $\B'$ and $\B''$ each determine a canonically metrized
holomorphic line bundle over $\M_g(L)$. These metrized line bundles
depend only on the variations $\B/\Z(1)$ and $W_{-1}\B$, and are
therefore isomorphic. Denote this common line bundle by
$\calB \to \M_g(L)$. The biextensions $\B'$ and $\B''$ determine
(and are determined by) meromorphic sections $s'$ and $s''$ of $\calB$,
respectively. There is therefore a meromorphic function $h$ on $\M_g(L)$
such that $s'' = hs'$. It follows from the main result of \cite{lear}
that this function is a rational function. (The philosophy is that
period maps of variations of mixed Hodge structure behave well at the
boundary.) The result follows as
$$
\nu(B''_C) = \log||s''(C)|| = \log|h(C)| + \log||s'(C)|| = \nu(B'_C)
+ \log |h(C)|.
$$

To conclude the proof, we now explain how to proceed when the
extension (\ref{putative}) is not split as a $\Q$ variation. Write
$\K = \bT \oplus \bT'$, where $\bT$ is the trivial submodule of $\K$
and $\bT'$ is its orthogonal complement. This is a splitting in the
category of $\Q$ variations by (\ref{decomp}). It also follows from
(\ref{decomp}) that the restriction of (\ref{putative}) to $\bT'$ is
split. Consequently, there is an inclusion of mixed Hodge structures
$\bT' \hookrightarrow \B(Z,W)$. As above, we may replace
$\B(Z,W)$ by the biextension $\B'=\B(Z,W)/\bT'$ after rescaling
lattices. This only changes the height by a non-zero rational number.
The weight graded $-1$ quotient of $\B'$ is the sum of at most one
copy of $\V(\lambda_3)$ and a trivial variation of weight $-1$.

Now suppose that $B_1$ and $B_2$ are two biextensions. We can
construct a new biextension $B_1 + B_2$ from them as follows:
Begin by taking their direct sum. Pull this back along the diagonal
inclusion
$$
\Z \hookrightarrow \Z\oplus \Z = Gr^W_0 \left(B_1\oplus B_2\right)
$$
to obtain a mixed Hodge structure $B$ whose weight $-2$ graded
quotient is
$$
Gr^W_{-2}\left(B_1 \oplus B_2\right) = \Z(1) \oplus \Z(1).
$$
Push this out along the addition map
$$
\Z(1) \oplus \Z(1) \to \Z(1)
$$
to obtain the sought after biextension $B_1\boxplus B_2$.  The following
result follows directly from \cite[(3.2.11)]{hain:heights}.

\begin{proposition}
The height of $B_1\boxplus B_2$ is the sum of the heights of $B_1$ and
$B_2$.
\end{proposition}

The biextension $\B'$ is easily seen to be the sum, in this sense, of
two biextensions. The first is constant with weight $-1$ quotient
equal to the trivial variation $\bT$ and the second is a variation with
weight $-1$ quotient equal to $\bH/\K$, which is either zero or one
copy of $\V(\lambda_3)$. Since the height of a constant biextension
is a constant, the result follows from the computation of the height
of a biextension with weight $-1$ quotient $\V(\lambda_3)$ above.

\section{Results for Abelian Varieties}
\label{abelian}

Denote the quotient of Siegel space $\h_g$ of rank $g$ by a finite index
subgroup $L$ of $Sp_g(\Z)$ by $\A_g(L)$. This is the moduli space of
abelian varieties with a level $L$ structure. In this section we state
results for $\A_g(L)$ analogous to those in Sections \ref{normal_functions}
and \ref{heights}. The proofs are similar, but much simpler, and are
left to the reader.

We call a representation of $Sp_g$ {\it even} if it has a symmetric
$Sp_g$-invariant inner product, and {\it odd} if it has a
skew symmetric $Sp_g$-invariant inner product. It follows from Schur's
Lemma that every irreducible representation of $Sp_g$ is either
even or odd. The even ones occur as polarized variations of Hodge
structure of even weight over each $\A_g(L)$, while the odd ones occur
as polarized
variations of Hodge structure only over $\A_g(L)$ of odd weight provided
$-I \notin L$. These facts are easily proved by adapting the arguments
in Section \ref{tech}.

The first theorem is the analogue of (\ref{fin_gen}) for abelian
varieties. It is similar to the result \cite{silverberg} of Silverberg.
The point in our approach is that $H^1(L,V)$ vanishes for all non-trivial
irreducible representations of $Sp_g$ by \cite{ragunathan}.

\begin{theorem}\label{norm_ab}
Suppose that $g\ge 2$ and that $L/\pm I$ is torsion free. If $\V \to \A_g(L)$
is a variation of Hodge structure of negative weight whose monodromy
representation is the restriction to $L$ of a rational representation of
$Sp_g$, then the group of generically defined normal functions associated
to this variation is finite.
\qed
\end{theorem}

Since there are no normal functions of infinite order over $\A_g(L)$,
we have the following
analogue of (\ref{height}). Suppose that $Z$ and $W$ are families of
homologically trivial cycles over $\A_g(L)$ in a family of smooth projective
varieties $p:X\to \A_g(L)$. Suppose that they are disjoint over the generic
point. Suppose further that $d + e = n-1$, where $d$, $e$ and $n$ are the
relative dimensions over $\A_g(L)$ of $Z$, $W$ and $X$,
respectively. Denote the fiber of $Z$ over $A\in \A_g(L)$ by $Z_A$, etc.

\begin{theorem}
If $g\ge 2$ and the monodromy of the local system $R^{2d+1}p_\ast \Q_X$
is the restriction to $L$ of a rational representation of $Sp_g$, then
there is a rational function $h$ on $\A_g(L)$ such that
$$
\langle Z_A, W_A \rangle = \log|h(A)|
$$
for all $A\in \A_g(L)$. \qed
\end{theorem}

One can formulate and prove analogues of these results for the moduli
spaces $\A_g^n(L)$ of abelian varieties of dimension $g$, $n$ marked
points, and a level $L$ structure.

We conclude this section with a discussion of Nori's results and their
relation to Theorems \ref{comp} and \ref{norm_ab}. We first recall
the main result of the last section of Nori's paper \cite{nori}.

\begin{theorem}[Nori]\label{nori_thm}
Suppose that $X$ is a variety that is an unbranched covering of a
Zariski open subset $U$ of $\A_g(L)$, where $L$ is torsion free.
Suppose that $\V$ is a variation of Hodge structure of negative weight
over $X$ that is pulled back from the canonical variation over $\A_g(L)$
of the same weight whose monodromy representation is irreducible and
has highest weight $\lambda$. Then the group of normal functions
defined on $X$ associated to this variation is finite unless
$$
\lambda = \begin{cases}
0 & \text{and $g \ge 2$ or;}\cr
\lambda_1 & \text{and $g\ge 3$ or;} \cr
\lambda_3 & \text{and $g=3$ or; }\cr
m_1\lambda_1 + m_2\lambda_2 & \text{$g=2$ and $m_1\ge 2$.}
\qed
\end{cases}
$$
\end{theorem}

This result may seem to contradict Theorem \ref{norm_ab}. The difference
can be accounted for by noting that Theorem \ref{norm_ab} only applies
to open subsets of the $\A_g(L)$, whereas Nori's theorem applies to
a much more general class of varieties which contains unramified
coverings of open subsets of the $\A_g(L)$. One instructive example is
$\M_3(l)$, where $l$ is odd and $\ge 3$. The map $\M_3(l) \to \A_3(l)$ is
branched along the hyperelliptic locus. Theorem \ref{norm_ab} does not
apply. However, Nori's Theorem \ref{nori_thm} does apply --- remember,
normal functions in weight $-1$ extend by (\ref{extend}). In this
way we realize the normal function associated to $\lambda_3$ in
Nori's result. Also,
by standard arguments, for each $n$, there is an open subset $U$
of $M_3(l)$ and an unbranched finite cover $V$ of $U$ over which the
natural projection $\M_3^n(l) \to \M_3(l)$ has a section. From this
one can construct $n$ linearly independent normal sections of the
jacobian bundle defined over $V$. Note that Nori's result does apply
to $V$, whereas (\ref{norm_ab}) does not.

\end{document}